\documentclass[english]{amsart}

\usepackage{enumerate}
\usepackage{amssymb}
\usepackage{array}
\usepackage{graphicx}
\usepackage{mathrsfs}
\usepackage{multicol}
\usepackage{mathtools}
\usepackage{ mathdots }
\usepackage{color}
\usepackage[usenames,dvipsnames,svgnames,table]{xcolor}
\usepackage{pgf,tikz,pgfplots}
\usepackage{tikz-cd}
\usetikzlibrary{arrows}
\pgfplotsset{compat=1.15}

\usepackage[curve,all]{xy}
\xyoption{matrix}
 \xyoption{curve}
 \xyoption{color}
 \xyoption{line}
\xyoption{arc}

\usepackage[shortlabels]{enumitem}
\setlist[enumerate]{leftmargin=7mm,topsep=0pt,itemsep=-1ex,partopsep=1ex,parsep=1ex,label=\rm{(\roman*)}}
\setlist[itemize]{leftmargin=5mm,topsep=0pt,itemsep=-1ex,partopsep=1ex,parsep=1ex,label=\raisebox{0.25ex}{\tiny$\bullet$}}

\usepackage[backref=false, colorlinks, linktocpage, citecolor = blue, linkcolor = blue]{hyperref}

\theoremstyle{plain}
\newtheorem{theorem}{Theorem}[section]

\newtheorem*{theoremaux}{Theorem \theoremauxnum}
\gdef\theoremauxnum{1}

\newtheorem*{main-theorem}{Main Theorem}
\newtheorem{proposition}[theorem]{Proposition}
\newtheorem*{propositionaux}{Proposition \propositionauxnum}
\gdef\propositionauxnum{1}

\newtheorem{lemma}[theorem]{Lemma}
\newtheorem*{lemmaaux}{Lemma \lemmaauxnum}
\gdef\lemmaauxnum{1}

\newtheorem{corollary}[theorem]{Corollary}

\newtheorem{question}{Question}
\newtheorem{problem}{Problem}
\newtheorem*{key-problem}{Key Problem}

\theoremstyle{definition}
\newtheorem{definition}[theorem]{Definition}

\newtheorem{example}[theorem]{Example}

\newtheorem*{works}{Related works}

\theoremstyle{remark}
\newtheorem{remark}[theorem]{Remark}

\makeatletter
\@namedef{subjclassname@2020}{%
  \textup{2020} Mathematics Subject Classification}
\makeatother

\newcommand{\leftexp}[2]{{\vphantom{#2}}^{#1}{#2}}

\newcommand{\incl}[1][r]{\ar@<-0.2pc>@{^(-}[#1] \ar@<+0.2pc>@{-}[#1]}

\newcommand{\hs}{\kern 0.8pt}

\renewcommand{\AA}{\mathcal{A}}
\newcommand{\BB}{\mathcal{B}}
\newcommand{\CC}{\mathcal{C}}
\newcommand{\AB}{\mathcal{AB}}

\newcommand{\Bl}{{\mathrm{Bl}}}

\renewcommand{\P}{\mathbb{P}}
\newcommand{\X}{\mathbb{X}}

\newcommand{\Hom}{\mathrm{Hom}}

\newcommand{\gr}{\mathrm{gr}}
\renewcommand{\div}{\mathrm{div}}

\newcommand{\W}{\mathcal{W}}

\newcommand{\Id}{\mathrm{Id}}
\newcommand{\Spec}{\mathrm{Spec}}
\newcommand{\Stab}{\mathrm{Stab}}
\newcommand{\inn}{\mathrm{inn}}

\newcommand{\V}{\mathcal{V}}

\renewcommand{\SS}{\mathcal{S}}

\newcommand{\Sone}{\mathbb{S}^1}
\renewcommand{\sl}{\mathfrak{sl}}

\newcommand{\g}{\mathfrak{g}}

\newcommand{\Q}{\mathbb{Q}}

\newcommand{\RR}{\mathcal{R}}

\newcommand{\Dyn}{\mathrm{Dyn}}
\newcommand{\Cone}{\mathrm{cone}}

\renewcommand{\H}{\mathrm{H}}

\newcommand{\K}{\mathrm{K}}

\newcommand{\A}{\mathbb{A}}
\newcommand{\iso}{\simeq}

\newcommand{\D}{\mathbb{D}}
\newcommand{\DDD}{\mathcal{D}}
\newcommand{\C}{\mathbb{C}}

\newcommand{\N}{\mathbb{N}}

\newcommand{\Z}{\mathbb{Z}}
\renewcommand{\K}{\mathbb{K}}
\newcommand{\R}{\mathbb{R}}
\newcommand{\G}{\mathbb{G}}
\newcommand{\Gm}{\mathbb{G}_m}

\renewcommand{\O}{\mathcal{O}}

\DeclareMathOperator{\SL}{SL}
\DeclareMathOperator{\SO}{SO}
\DeclareMathOperator{\SU}{SU}

\DeclareMathOperator{\Aut}{Aut}

\DeclareMathOperator{\PGL}{PGL}

\DeclareMathOperator{\Gal}{Gal}

\def\ga{\,^\gamma\hskip-1pt}

\title[Real structures on varieties with a reductive group action]{Real structures on complex algebraic varieties with a reductive group action}

\author[Ronan Terpereau]{Ronan Terpereau}
\thanks{The author is supported by the ANR Project FIBALGA ANR-18-CE40-0003-01.
This work received partial support from the French ``Investissements d\textquoteright Avenir" program and from project ISITE-BFC (contract ANR-lS-IDEX-OOOB). The IMB receives support from  the EIPHI Graduate School (contract ANR-17-EURE-0002).
}

\address{Institut de Math\'{e}matiques de Bourgogne, UMR 5584 CNRS, Universit\'{e} de Bourgogne, F-21000 Dijon, France}
\email{ronan.terpereau@u-bourgogne.fr}


\begin{document}

\begin{abstract}
This is a survey article devoted to the study of real structures on complex algebraic varieties endowed with a reductive group action.
\end{abstract}

\maketitle

\tableofcontents

\section*{Introduction}
All notions appearing in the introduction will be carefully defined in the later sections.
Let $X$ be a complex algebraic variety. 
Recall that a real algebraic variety $W$, together with an isomorphism $X \simeq W_{\C} := W \times_{\Spec(\R)} \Spec(\C)$ of complex algebraic varieties, is said to be a \emph{real form} of $X$. (Note that for some authors, the definition of real form is slightly different from ours and applies to $W$ instead of $X$; this does not change the Key Problem we will consider, it modifies only its formulation.) 
It is a classical problem in algebraic geometry to determine the real forms of a given complex algebraic variety which usually splits 	into two parts:
\begin{problem}[real forms in the classical setting]\label{prob1} \ 
\begin{enumerate}
\item Does the complex algebraic variety $X$ admit a real form?
\item If ``yes", how to classify them? Are there a finite number (up to isomorphism)?
\end{enumerate}
\end{problem}

However, in the present article we are mainly interested in algebraic varieties endowed with a reductive group action, and so only real forms with symmetries coming from the group action will be taken into account. Let us explain what this means. 
Let $G$ be a complex algebraic group, let $F$ be a real algebraic group that is a real form of $G$ (in the category of algebraic groups), and let now $X$ be a complex \emph{$G$-variety} (i.e.~an algebraic variety endowed with a faithful $G$-action). A real algebraic $F$-variety $W$ together with an isomorphism $X \simeq W_{\C}$ of complex algebraic $G$-varieties is said to be a real $F$-form of $X$.
In this setting, Problem \ref{prob1} admits the following reformulation.
\begin{problem}[real forms in the equivariant setting]\label{prob2}  \ 
\begin{enumerate}
\item Does the complex $G$-variety $X$ admit a real $F$-form?
\item If ``yes", how to classify them? Are there a finite number (up to $F$-equivariant isomorphism)?
\end{enumerate}
\end{problem} 

It follows from the work of Weil \cite{Wei56} and Borel--Serre \cite{BS64} that there is a correspondence between the notions of real forms and of \emph{real structures} (i.e.~antiregular involutions). More precisely, real forms of complex algebraic groups correspond to \emph{real group structures} (see Definition \ref{def:real group structure}) and
real forms of complex algebraic varieties endowed with an algebraic group action correspond to \emph{effective equivariant real structures} (see Definition \ref{def: equiv real structure}). This correspondence is usually referred to as \emph{Galois descent}. 
In practice, it is more convenient to work with real structures instead of real forms; indeed, real structures can be extended or restricted with respect to a dense open subset, and they are parametrized by some Galois cohomology pointed set (see Proposition \ref{prop:Galois H1 to param eq real structures}). 
Consequently, when studying real forms of complex algebraic varieties endowed with an algebraic group action, we are naturally led to consider the following Key Problem instead of Problem \ref{prob2}.

\begin{key-problem}[real structures in the equivariant setting] \ 
\begin{enumerate}
\item\label{item: existence part of Key-Prob3} Given a real group structure $\sigma$ on $G$, does there exist an effective $(G,\sigma)$-equivariant real structure on the complex $G$-variety $X$? 
\item\label{item: quantity part of Key-Prob3} If ``yes", how to classify them? Are there a finite number (up to equivalence)?
\end{enumerate}
\end{key-problem} 
 
The goal of this survey article is to succinctly present the current state of the art concerning the Key Problem, with a special focus on progress made over the past decade. Concretely, we consider the Key Problem for several families of $G$-varieties \footnote{Throughout this article, varieties are assumed to be normal as this assumption is crucial in the different combinatorial descriptions that will appear.} for which there exist well-established combinatorial descriptions. Since each of these descriptions is however quite technical, we do not to go into details, but rather try to focus on the big picture.
We will see that each time, a complete answer to Question \ref{item: existence part of Key-Prob3} of the Key Problem can be given in terms of certain stability properties satisfied by the combinatorial data on which the Galois group $\Gal(\C/\R)$ acts naturally (sometimes with an additional cohomological condition), but that in most cases Proposition \ref{prop:Galois H1 to param eq real structures} is the best answer that can be made to Question \ref{item: quantity part of Key-Prob3} of the Key Problem.

\subsection*{Structure of the article}
Section \ref{sec: real group structures on reductive groups} contains preliminaries on real group structures on reductive groups and a definition of the famous $\star$-action of the Galois group $\Gal(\C/\R)$ on the based root datum associated to a triple $(G,B,T)$ with $B \subseteq G$ a Borel subgroup and $T \subseteq B$ a maximal torus (see Definition \ref{def: star-action}); this action will play a crucial role when answering Question \ref{item: existence part of Key-Prob3} of the Key Problem.
Section \ref{sec: equivariant real structures} contains preliminaries on equivariant real structures on $G$-varieties, with particular emphasis on the case of almost homogeneous $G$-varieties.
In Section \ref{sec: warm-up}, we then consider the Key Problem for two classical families of complex algebraic varieties with a reductive group action, namely the linear representations and the nilpotent orbits (and the normalizations of their closures) in semisimple Lie algebras. 
Section \ref{sec: real structures on T-varieties} is devoted to varieties with a torus action; we first consider the case of toric varieties, then the case of arbitrary affine $T$-varieties.
In Section \ref{sec: real structures on spherical varieties}, we pass from the toric case to the spherical case; this generalization is done in two steps: we first consider spherical homogeneous spaces, then their equivariant embeddings.
In Section \ref{sec:almost homg SL2-threefolds}, we consider the Key Problem for the almost homogeneous $\SL_2$-threefolds for which, contrary to the spherical case, the $\star$-action is insufficient to give a combinatorial solution to Question \ref{item: existence part of Key-Prob3} of the Key Problem.
Finally, Section \ref{sec: open questions} contains some open questions and lines of research related to  the Key Problem.

\subsection*{Foreword}
While I have endeavored to be comprehensive, the decisions made to keep this article concise inevitably reflect my personal interests.

\subsection*{Acknowledgments}
I would like to thank Michel Brion, Adrien Dubouloz, Pierre-Alexandre Gillard and Lucy Moser-Jauslin for their valuable suggestions on a former version of this article. I also thank Thibaut Delcroix for providing me with the homogeneous spherical data for Examples \ref{example A} and \ref{example B}.
I am very grateful to the referee for the careful reading of the paper and for his/her list of recommendations.

\subsection*{Notation}
In this article, a \emph{variety} (over a field $\K \in \{ \R,\C\}$) is a separated scheme of finite type (over $\K$) that is geometrically integral  and \textbf{normal}, and an \emph{algebraic group} (over a field $\K$) is a finite type group scheme (over $\K$). By an \emph{algebraic subgroup}, we always mean a closed subgroup scheme. 
A connected linear algebraic group $G$ is called \emph{reductive} if its largest connected unipotent normal subgroup is trivial. We refer to \cite{Hum75} for a general reference on linear algebraic groups, and to \cite{Con14} for a more advanced reference on reductive group schemes.

We denote by $G$ a complex reductive group, and we will always assume (except in Sections \ref{sec: real group structures on reductive groups} and \ref{sec: equivariant real structures}) that $G$ is of \emph{simply-connected type}, i.e.~isomorphic to the product of a torus and a simply-connected semisimple group.
We denote by $Z(G)$ the center of $G$, by $B$ a Borel subgroup of $G$, by $T$ a maximal torus of $G$ contained in $B$, and by $U$ the unipotent radical of $B$ (which is also a maximal unipotent subgroup of $G$). 
We write $\X=\X(T)=\Hom_{\gr}(T,\Gm)$ for the character lattice of $T$, $\X^\vee=\X^\vee(T)=\Hom_{\gr}(\Gm,T)$ for the cocharacter lattice of $T$, and $\X_\Q:=\X \otimes_\Z \Q$ and  $\X_\Q^\vee:=\X^\vee \otimes_\Z \Q$ for the corresponding $\Q$-vector spaces.
If $H$ is an algebraic subgroup of $G$, then $N_G(H)$ denotes the  normalizer of $H$ in $G$. 
When $G$ is semisimple, we denote its Dynkin diagram by $\Dyn(G)$.

For $c \in G$, we write \[\inn_c\colon G \to G,\ g \mapsto cgc^{-1}.\] We say that $\inn_c$ is an \emph{inner} automorphism of $G$. An automorphism of $G$ that is not inner is called \emph{outer}.
The set of inner automorphisms of $G$ is a normal subgroup of $\Aut_{\gr}(G)$ whose quotient group identifies with $\Aut(\Dyn(G))$ when $G$ is semisimple of simply-connected type.

The Galois group of the extension $\C/\R$ is denoted by 
\[
\Gamma:=\Gal(\C/\R)=\{\Id,\gamma\} \simeq \Z/2\Z.
\]
If $A$ is a $\Gamma$-group, then the first Galois cohomology pointed set associated to $A$ is denoted by $\H^1(\Gamma,A)$.  
If furthermore  $A$ is an abelian group, then we can also consider $\H^2(\Gamma,A)$, and both $\H^1(\Gamma, A)$ and $\H^2(\Gamma,A)$ are abelian groups. 
We refer to \cite{Ser02} for a general reference on Galois cohomology (or to \cite[Section 1]{Gil23} for a detailed summary on Galois descent), and to \cite{Man20} for a general reference on real algebraic geometry.

\section{Real group structures on reductive groups}\label{sec: real group structures on reductive groups}
In this section, we introduce the notion of \emph{real group structures} on complex algebraic groups. In particular, since we are mainly interested in the case of reductive groups, we explain how to obtain all real group structures on complex reductive groups by piecing together real group structures on complex tori and on complex simply-connected simple groups. Then we recall how the choice of a real group structure on a reductive group induces a $\Gamma$-action (the so-called \emph{$\star$-action}) on its based root datum; this $\star$-action will play an important role when looking at the Key Problem in the next sections for linear representations, nilpotent orbits in semisimple Lie algebras, $T$-varieties, and spherical varieties.

\begin{definition}Let $G$ be a complex algebraic group.  \label{def:real group structure}
\begin{enumerate}
\item A \emph{real form} of $G$ is a pair $(F,\Theta)$ with $F$ a real algebraic group and $\Theta\colon\  G \to F_\C$ an isomorphism of complex algebraic groups. (Most of the time we will drop the isomorphism $\Theta$ and simply write that $F$ is a form of $G$, but it is important to remember that a real form of $G$ is a pair, and not just a real algebraic group.)
\item A \emph{real group structure} $\sigma$ on $G$ is a scheme involution on $G$ such that the diagram 
\[
\xymatrix@R=4mm@C=2cm{
    G \ar[rr]^{\sigma} \ar[d]  && G \ar[d] \\
    \Spec(\C)  \ar[rr]^{\Spec(z \mapsto \overline{z})} && \Spec(\C)  
  }
  \]
commutes and 
 \[\iota_G \circ \sigma=\sigma \circ \iota_G \ \ \text{ and }\ \ m_G \circ (\sigma \times \sigma)=\sigma \circ m_G,\] where $\iota_G\colon G \to G$ is the inverse morphism and $m_G\colon G \times G \to G$ is the multiplication morphism.
\item Two real group structures $\sigma$ and $\sigma'$ on $G$ are \emph{equivalent} if there exists a complex algebraic group automorphism $\psi \in \Aut_{\gr}(G)$ such that $\sigma'=\psi \circ \sigma \circ \psi^{-1}$.
They are \emph{strongly equivalent} if we can choose $\psi=\inn_c$ for some $c \in G$.
\item If $\sigma_1$ and $\sigma_2$ are two real group structures on $G$, we may write $\sigma_2=\varphi \circ \sigma_1$ with $\varphi \in \Aut_{\gr}(G)$, and we say that $\sigma_2$ is an \emph{inner} (resp.~an \emph{outer}) \emph{twist} of $\sigma_1$ if $\varphi$ is an inner (resp.~an outer) automorphism of $G$. Let us note that the relation ($\sigma_2$ is an inner twist of $\sigma_1$) is an equivalence relation on the set of real group structures on $G$.
\item If $G$ is a complex algebraic group with a real group structure $\sigma$, then the set of fixed points $G(\C)^\sigma$ is called the \emph{real locus} (or \emph{real part}) of $(G,\sigma)$.
\end{enumerate}
\end{definition}

\begin{remark}
The notions of real forms and real group structures are particular cases of the notions of $\mathbb{K}$-forms and descent data (in the category of algebraic groups) in the case where the Galois extension considered $\mathbb{K}'/\mathbb{K}$ is $\C/\R$. The interested reader is referred to \cite[Section 2.1]{MJT20} to learn more about these notions.
\end{remark}

There is a correspondence between real group structures on $G$ and real forms of $G$ given as follows (see also \cite[Section 1.4]{FSS98}).
\begin{itemize}
\item 
If $F$ is a real form of $G$, then the composition $G \simeq F_\C  \xrightarrow{\Id \times \Spec(z \mapsto \overline{z})}  F_\C \simeq G$ is a real group structure on $G$.
\item Conversely, if $\sigma$ is a real group structure on $G$, then the categorical quotient $F:=G/\langle \sigma \rangle$ gives a real form of $G$; an isomorphism $G \simeq F_\C $  being induced by $(q,f)$, where $q\colon\  G\to F$ is the quotient morphism and $f\colon\  G \to \Spec(\C)$ is the structure morphism.
\item Moreover, two real forms of $G$ are isomorphic (as real algebraic groups) if and only if the corresponding real group structures are equivalent. 
\end{itemize}

\begin{remark}
The real locus of $(G,\sigma)$ identifies with the group of real points of the real algebraic group $G/\langle \sigma \rangle$, whence the name. Moreover, since the group of real points of a real algebraic group inherits a natural structure of a real Lie group, the real locus $G(\C)^\sigma$ is a real Lie group.
\end{remark}

\smallskip

We now assume that $G$ is a complex reductive group. There exists a central isogeny $\widetilde{G} \to G$, where $\widetilde{G}$ is the direct product of a torus $T$ and of a simply-connected semisimple group $G'$, such that every real group structure $\sigma$ on $G$ lifts uniquely to a real group structure $\widetilde{\sigma}$ on $\widetilde{G}$ (see \cite[Section 1.1]{MJT18} for details), and $\widetilde{\sigma}$ stabilizes the two factors $T$ and $G'$. Hence, determining real group structures on complex reductive groups reduces to determining real group structures on complex tori and on complex simply-connected semisimple groups.

Both classifications are well-known: for complex tori it is recalled in Lemma \ref{lem: real form on tori} below, and for simply-connected semisimple groups, it reduces to the case of simply-connected simple groups (see Lemma \ref{lem:easy_lemma_reduction}), in which case real group structures are classified in terms of Dynkin diagrams enriched by blackening a subset of the nodes and connecting some pairs of vertices by arrows (the \emph{Satake diagrams} or the \emph{Vogan diagrams}); see e.g.~\cite[Section 1.7.2]{GW09} or \cite[Section V\!I.10]{Kna02} for a summary of the classification of real group structures on simple groups.   

\begin{lemma} \label{lem: real form on tori} \emph{(see \cite[Lemma 1.5]{MJT18} and \cite[Theorem 2]{Cas08})}\\
Let $T \simeq \G_{m,\C}^n$ be an $n$-dimensional complex torus.
\begin{enumerate}
\item If $n=1$, then $T$ has two inequivalent real group structures, defined by $\sigma_0\colon t \mapsto \overline{t}$, corresponding to $\G_{m,\R}$, and by $\sigma_1\colon t \mapsto \overline{t^{-1}}$, corresponding to the circle $\Sone$. 
\item If $n=2$, then $\sigma_2\colon  (t_1,t_2) \mapsto (\overline{t_2},\overline{t_1})$ defines a real group structure on $T$ corresponding to the Weil restriction $R_{\C/\R}(\G_{m,\C})$ (we refer to \cite[Section 7.6]{BLR90} for a detailed treatment of the Weil restriction functor).
\item \label{item: n at least 2} If $n\ge 1$,  then every real group structure on $T$ is equivalent to exactly one real group structure of the form  $\sigma_0^{\times n_0}\times\sigma_1^{\times n_1}\times\sigma_2^{\times n_2}$, where $n=n_0+n_1+2n_2$. In other words, every real torus is isomorphic to $\G_{m,\R}^{n_0} \times (\Sone)^{n_1} \times R_{\C/\R}(\G_{m,\C})^{n_2}$.
\end{enumerate}
\end{lemma}

\begin{lemma} \emph{(see \cite[Lemma 1.7]{MJT18})}\label{lem:easy_lemma_reduction}
Let $\sigma$ be a real group structure on a complex simply-connected semisimple group $G' \simeq \prod_{i \in I} G_i$, where the $G_i$ are the simple factors of $G'$. Then, for a given $i \in I$, we have the following possibilities:
\begin{enumerate}[leftmargin=4mm]
\item $\sigma(G_i)=G_i$ and $\sigma_{|G_i}$ is a real group structure on $G_i$; or
\item there exists $j \neq i$ such that $\sigma(G_i)=G_j$, then $G_i \simeq G_j$ and $\sigma_{| G_i \times G_j}$ is equivalent to $(g_1,g_2) \mapsto (\sigma_0(g_2),\sigma_0(g_1))$ for a real group structure $\sigma_0$ on $G_i \simeq G_j$.
\end{enumerate}
\end{lemma}

\begin{example}\label{ex:real group structures on SL2}
There are two inequivalent real group structures on $\SL_2$ given by $\sigma_s(g)=\overline{g}$, whose real locus is $\SL_2(\R)$, and its inner twist $\sigma_c(g)={}^t{\overline{g}^{-1}}$, whose real locus is $\SU_2(\C)$. 
(Here $\overline{g}$ denotes the complex conjugate.)
There are four inequivalent real group structures on $\SL_2 \times \SL_2$ given by $\sigma_i \times \sigma_j$, with $(i,j) \in \{(s,s),(s,c),(c,c)\}$, and $\sigma\colon  (g_1,g_2) \mapsto (\sigma_s(g_2),\sigma_s(g_1))$, which is equivalent to the real group structure $(g_1,g_2) \mapsto (\sigma_c(g_2),\sigma_c(g_1))$ and whose real locus is $\SL_2(\C)$.
Similarly, we let the reader check that there are six inequivalent real group structures on $\SL_2 \times \SL_2 \times \SL_2$ and nine inequivalent real group structures on $\SL_2 \times \SL_2 \times \SL_2 \times \SL_2$.
\end{example}

\smallskip

Let $B \subseteq G$ be a Borel subgroup, and let $T \subseteq B$ be a maximal torus. We say that $(T,B)$ is a \emph{Borel pair} in $G$.
Let $\sigma$ be a real group structure on $G$.
Then $(T',B'):=(\sigma(T),\sigma(B))$ is another Borel pair in $G$, and so there exists $c_{\sigma} \in G$ (unique up to left multiplication by an element of $T$) such that 
\[
c_{\sigma} T' c_{\sigma}^{-1}=T\ \ \ \ \text{and}\ \ \ \ c_{\sigma} B' c_{\sigma}^{-1}=B.
\]
Let $\theta:=\inn_{c_{\sigma}} \circ \sigma$ that is an antiregular automorphism of $G$. 
Since $\theta(T)=T$ and $\theta(B)=B$, the automorphism $\theta$ induces lattice automorphisms of $\X=\X(T)$ and $\X^\vee=\X^\vee(T)$ defined as follows:
 \[
 \forall \chi \in \X,\ \chi \mapsto\ga \chi :=\sigma_0 \circ \chi \circ \theta^{-1}\ \ \ \ \text{ and } \ \ \ \ \forall \lambda \in \X^\vee,\  \lambda \mapsto \ga \lambda:=\theta  \circ \lambda \circ \sigma_0,\]
where $\sigma_0(t)=\overline{t}$ is the complex conjugation on $\G_{m,\C}$. One can check that these two lattice automorphisms are of order $\leq 2$, that they do not depend on the choice of $c_{\sigma}$, and that they stabilize the sets of roots $\RR \subseteq \X$, simple roots $\SS \subseteq \RR$, coroots $\RR^\vee \subseteq \X^\vee$, and simple coroots $\SS^\vee \subseteq \RR^\vee$ associated with the triple $(G,B,T)$; see e.g.~\cite[Remark 7.1.2]{Con14} or \cite[Section A.2]{MJT18} and references therein for details.

\begin{definition}\label{def: star-action}
With the previous notation, the $\Gamma$-action on the \emph{based root datum} $(\X,\X^\vee,\RR,\RR^\vee,\SS,\SS^\vee)$ associated with $(G,B,T)$ is called the \emph{$\star$-action} induced by $\sigma$; it does not depend on the choice of the Borel pair $(T,B)$ in $G$.
\end{definition} 

 \begin{example}\label{ex: star-action for tori}
If $G=B=T$ and $\sigma$ is a real group structure on $T$, then the corresponding $\star$-action is given by
 \[
 \forall \chi \in \X,\ \chi \mapsto\ga \chi =\sigma_0 \circ \chi \circ\sigma \ \ \ \ \text{ and } \ \ \ \ \forall \lambda \in \X^\vee,\  \lambda \mapsto \ga \lambda=\sigma  \circ \lambda \circ \sigma_0.\]
 \end{example}

\begin{remark}\label{rk: star action trivial in the split case}
It follows from the definition of the $\star$-action that if $\sigma$ is an inner twist of $\sigma'$, then the $\star$-actions induced by $\sigma$ and $\sigma'$ coincide.
In particular, if $\sigma$ is an inner twist of a \emph{split} real group structure $\sigma'$ on $G$ (i.e.~$G$ contains a maximal torus $T_0$ preserved by $\sigma'$ and such that $\sigma'_{|T_0}=\sigma_0^{\times \dim(T_0)}$, with the notation of Lemma \ref{lem: real form on tori}), then the $\star$-action induced by $\sigma$ is trivial.
On the other hand, if $\sigma$ is an outer twist of a split real group structure on $G$, then it follows from \cite[Section 1.5, second half of p.~41]{Con14} that the $\star$-action induced by $\sigma$ is nontrivial.
\end{remark}

\section{Equivariant real structures on \texorpdfstring{$G$}{G}-varieties}\label{sec: equivariant real structures}
In this section, we first introduce the notion of \emph{equivariant real structures} on algebraic varieties endowed with an algebraic group action, then we give a general criterion for the existence of an equivariant real structure in the homogeneous case, and finally we recall how Galois cohomology can be used to parametrize equivalence classes of equivariant real structures. 

\smallskip

Let $G$ be a complex algebraic group, let $F$ be a real form of $G$, and let $\sigma$ be the corresponding real group structure on $G$.

\begin{definition}\label{def: equiv real structure}
Let $X$ be a complex algebraic $G$-variety.
\begin{itemize}
\item A \emph{real $F$-form} of $X$ is a pair $(W,\Xi)$ with $W$ a real algebraic $F$-variety and $\Xi\colon X \to W_\C$ an isomorphism of complex algebraic $G$-varieties. (As for real forms of complex algebraic groups, we will most of the time drop the isomorphism $\Xi$ and simply write that $W$ is a form of $X$.) 
\item A \emph{$(G,\sigma)$-equivariant real structure} on $X$ is an antiregular involution $\mu$ on $X$, that is, a scheme involution on $X$  such that the following diagram commutes
\[
\xymatrix@R=4mm@C=2cm{
    X \ar[rr]^{\mu} \ar[d]  && X \ar[d] \\
    \Spec(\C)  \ar[rr]^{\Spec(z \mapsto \overline{z})} && \Spec(\C)  
  }
  \]
and that satisfies the condition
\begin{equation}\label{eq:def eq real structure}
\forall g \in G, \ \forall x \in X,\ \ \mu(g \cdot x)=\sigma(g) \cdot \mu(x). 
\end{equation}
\item A $(G,\sigma)$-equivariant real structure $\mu$ on $X$ is \emph{effective} if $X$ is covered by $\Gamma$-stable affine open subsets ($\Gamma$ acting on $X$ through $\mu$), i.e.~if the categorical quotient $X/\langle \mu \rangle$, which always exists as a real algebraic space, is in fact a real algebraic variety (see \cite[Proposition~V.1.8]{SGA1}). 
\item Two $(G,\sigma)$-equivariant real structures $\mu$ and $\mu'$ on $X$ are \emph{equivalent} if there exists a $G$-equivariant automorphism $\varphi \in \Aut^{G}(X)$ such that $\mu'=\varphi \circ \mu\circ \varphi^{-1}$. 
\item The \emph{real locus} (or \emph{real part}) of a $(G,\sigma)$-equivariant real structure $\mu$ on $X$ is the (possibly empty) set  of fixed points $X(\C)^\mu$; it identifies with the set of $\R$-points of the corresponding real form $X/\left \langle \mu \right \rangle$ and is furthermore endowed with an action of the real Lie group $G(\C)^\sigma$.
\end{itemize}
\end{definition}

\begin{remark}
The notions of real forms and equivariant real structures are particular cases of the notions of $\mathbb{K}$-forms and descent data (in the category of algebraic varieties endowed with an algebraic group action) in the case where the Galois extension considered $\mathbb{K}'/\mathbb{K}$ is $\C/\R$. The interested reader is referred to \cite[Section 2.2]{MJT20} to learn more about these notions.
\end{remark}

As for algebraic groups (see Section \ref{sec: real group structures on reductive groups}), there is a bijection between isomorphism classes of real $F$-forms of $X$ (as real algebraic $F$-varieties) and equivalence classes of effective $(G,\sigma)$-equivariant real structures on $X$ (see \cite[Section 5]{Bor20}). 

\begin{remark}
Assume that the complex algebraic variety $X$ has a real structure; it induces a $\Gamma$-action on $X$.
If $X$ is quasiprojective or covered by $\Gamma$-stable quasiprojective open subsets, then $X$ is covered by $\Gamma$-stable affine open subsets. 
In particular, since homogeneous spaces under the action of a linear algebraic group are quasiprojective (this follows from a theorem of Chevalley; see \cite[Section 11.2, Theorem]{Hum75}), equivariant real structures on homogeneous spaces under a linear algebraic group are always effective.
\end{remark}

The following lemma gives a general criterion for the existence of an equivariant real structure in the case where the algebraic $G$-variety $X$ is a homogeneous space. 
(This criterion is however not so useful in real life as condition \ref{eq: involution} of Lemma \ref{lem: two conditions} is generally difficult to check.)

\begin{lemma} \emph{(\cite[Lemma 2.4]{MJT18})}\label{lem: two conditions}
Let  $X=G/H$ be a homogeneous space.
Then $X$ has a $(G,\sigma)$-equivariant real structure if and only if there exists $t \in G$ such that
\begin{enumerate}
\item \label{eq: sigma compatible}
 $\sigma(H)=tHt^{-1}$ \ $((G,\sigma)$\emph{-compatibility condition}); and
\item \label{eq: involution}
$\sigma(t)t \in H$\ \ (\emph{involution condition});
\end{enumerate}
in which case a $(G,\sigma)$-equivariant real structure on $X$ is given by 
\[\forall k\in G,\ \mu(kH):=\sigma(k)tH.\]
\end{lemma}

\begin{example} \label{ex:G/P 1}
Let $G$ be a complex reductive group with a real group structure $\sigma$, and let $X=G/P$ be a flag variety. Then Lemma \ref{lem: two conditions} implies that $X$ has a $(G,\sigma)$-equivariant real structure if and only if $\sigma(P)$ is conjugate to $P$. This last condition can be interpreted combinatorially. Let $B \subseteq G$ be a Borel subgroup and let $T \subseteq B$ be a maximal torus, and consider the $\star$-action induced by $\sigma$ on the based root datum $(\X,\X^\vee,\RR,\RR^\vee,\SS,\SS^\vee)$ associated to the triple $(G,B,T)$.
Conjugacy classes of parabolic subgroups of $G$ are in bijection with the powerset of $\SS$, and one can check that $\sigma(P)$ is conjugate to $P$ if and only if the subset of simple roots corresponding to the conjugacy class of $P$ is $\Gamma$-stable.

Moreover, if a $(G,\sigma)$-equivariant real structure exists on $X$, then it is unique up to equivalence (this will be a consequence of Proposition \ref{prop:Galois H1 to param eq real structures} since the group of $G$-equivariant automorphisms $\Aut^G(X)$ is trivial when $X$ is a flag variety). 

Let us consider in more detail the case $G=\SL_3$. 
Let $c=c^{-1}=\begin{bmatrix}
0 & 0 & -i\\ 0 & -1 & 0 \\ i & 0 & 0
\end{bmatrix}\in \SL_3$.
Up to equivalence, the group $\SL_3$ has three real group structures given by $\sigma_0(g)=\overline{g}$, which is split and whose real locus is $\SL_3(\R)$,
$\sigma_1(g)=c(\leftexp{t}{\overline{g}}^{-1})c^{-1}$, which is \emph{quasi-split} (i.e.~stabilizes a Borel subgroup of $G$) and whose real locus is $\SU(1,2)$, and  $\sigma_2(g)=\leftexp{t}{\overline{g}}^{-1}=\inn_c \circ \sigma_1$, which is an inner twist of $\sigma_1$ and whose real locus is $\SU(3)$.
Let $\SS=\{\alpha_1,\alpha_2\}$ be the set of simple roots of $G$. Let $P_{\alpha_1}$ and $P_{\alpha_2}$ be representatives of the the conjugacy classes of parabolic subgroups of $G$ associated with the singletons $\{\alpha_1\}$ and $\{\alpha_2\}$ respectively. On the one hand, if $\sigma$ is equivalent to $\sigma_0$, then the corresponding $\star$-action is trivial, and so 
$G/P_{\alpha_i}$ admits a unique, up to equivalence, $(G,\sigma)$-equivariant real structure. On the other other hand, if $\sigma$ is equivalent to $\sigma_1$ or $\sigma_2$, then  the corresponding $\star$-action swaps $\alpha_1$ and $\alpha_2$, and so $G/P_{\alpha_i}$ has no $(G,\sigma)$-equivariant real structure.
\end{example}

Let us mention that if $\sigma$ and $\sigma'$ are two real group structures on $G$ conjugate by an inner automorphism, then there is a bijection between the equivalence classes of $(G,\sigma)$-equivariant real structures and of $(G,\sigma')$-equivariant real structures on $X=G/H$ (see \cite[Proposition 1.1]{MJT20}). On the other hand, if $\sigma$ and $\sigma'$ are conjugate by an outer automorphism, then these two sets may have different cardinalities.

\begin{example}\label{ex:counter-ex not strongly eq}
Let $G=\mathbb{G}_{m,\C}^{2}$, let $H=\{1\} \times \mathbb{G}_{m,\C}$, and let 
\[
\begin{array}{rcll}
\sigma\ \colon & G \to G,&\ (u,v) \mapsto (\overline{u},\overline{v^{-1}});&\\ 
\varphi\ \colon & G \to G,&\ (u,v) \mapsto (uv,v);&\ \text{and}\\ 
\sigma'=\varphi \circ \sigma \circ \varphi^{-1}\ \colon & G \to G, &\ (u,v) \mapsto (\overline{u} \overline{v^{-2}},\overline{v^{-1}}).&
\vspace{-1mm}
\end{array}
\]
Then $\sigma(H)=H$ but $\sigma'(H)=\{(t^2,t)\ | \ t \in \mathbb{G}_{m,\C}\} \neq H$, and so Lemma \ref{lem: two conditions} implies (because $G$ is commutative) that $X=G/H$ admits a $(G,\sigma)$-equivariant real structure but no $(G,\sigma')$-equivariant real structure. 
\end{example}

If we now suppose that $X$ is an \emph{almost homogeneous} $G$-variety, i.e.~ an algebraic $G$-variety with a dense open orbit $X_0=G/H$, then condition \eqref{eq:def eq real structure} implies that an equivariant real structure $\mu$ on $X$ stabilizes $X_0$, and so $\mu_{|X_0}$ is an equivariant real structure on $X_0$. 
But an equivariant real structure on $X_0$ need not extend to $X$ (however, if such an extension exists, then it is unique).
That is why, when studying equivariant real structures on almost homogeneous varieties, we are naturally led to first study equivariant real structures on homogeneous spaces. 
Let us note that the restriction map $\mu \mapsto \mu_{|X_0}$ induces a map
\begin{align*}
\Upsilon \colon &\{\text{equivalence classes of $(G,\sigma)$-equivariant real structure on $X$}\} \\
&\to  \{\text{equivalence classes of $(G,\sigma)$-equivariant real structure on $X_0$}\}
\end{align*}
that is generally neither surjective nor injective (see Example \ref{ex:first bis example}).
Indeed, $\Upsilon$ fails to be surjective when $X_0$ admits a $(G,\sigma)$-equivariant real structure $\mu$ that does not extend (nor any $(G,\sigma)$-equivariant real structure equivalent to $\mu$) to $X$, and $\Upsilon$ fails to be injective when inequivalent $(G,\sigma)$-equivariant real structures $\mu$ and $\mu'$ on $X$ become equivalent after restriction to $X_0$, i.e.~ when any $\varphi \in \Aut^G(X_0)$ such that $\mu'_{|X_0}=\varphi \circ \mu_{|X_0} \circ \varphi^{-1}$ does not extend to an element of $\Aut^G(X)$.
The map $\Upsilon$ is however injective when the natural injective homomorphism  $\Aut^{G}(X) \hookrightarrow \Aut^{G}(X_0)$ is an isomorphism (which is, for instance, the case for toric varieties).

Let us mention that a general strategy to determine the equivariant real structures on almost homogeneous varieties, relying on Luna-Vust theory, is detailed in \cite[Section 3.5]{MJT20}. It is this strategy that will be applied to determine the equivariant real structures on toric varieties (in Section \ref{subsec: real structures on toric varieties}), on spherical embeddings (in Section \ref{subsec: extension of real structures to sph embeddings}), and on almost homogeneous $\SL_2$-threefolds (in Section \ref{sec:almost homg SL2-threefolds}).

\smallskip

Once we know the existence of a $(G,\sigma)$-equivariant real structure $\mu$ on an algebraic $G$-variety $X$, Galois cohomology can be used to parametrize the equivalence classes of all of them. 
First, observe that the group of $G$-equivariant automorphisms $\Aut^{G}(X)$ is endowed with a $\Gamma$-group structure as follows:
\[
 \Gamma \times \Aut^{G}(X) \to \Aut^{G}(X) ,\ (\gamma,\varphi) \mapsto
  \mu \circ \varphi \circ \mu.
\]
The next result is then a straightforward consequence of the definition of cohomologous $1$-cocycles with values in the $\Gamma$-group $\Aut^{G}(X)$.

\begin{proposition} \label{prop:Galois H1 to param eq real structures} \emph{(see \cite[Corollary 8.2]{Wed18} or \cite[Lemma 2.11]{MJT18})}\\
Let  $X$ be a complex algebraic $G$-variety with a $(G,\sigma)$-equivariant real structure $\mu_0$. The Galois group $\Gamma$ acts on $\Aut^G(X)$ by $\mu_0$-conjugacy as above. Then the map
\[\begin{array}{ccl}
\hspace{-4mm} \H^1(\Gamma,\Aut^G(X)) &\to &\{ \text{equivalence classes of $(G,\sigma)$-equivariant real structures on $X$}\}\\
 \varphi &\mapsto &\ \ \ \varphi \circ \mu_0
 \end{array}\]
 is a bijection that sends the trivial cohomology class to the equivalence class of $\mu_0$.
\end{proposition}

\begin{corollary}\emph{(see \cite[Corollary 2.8]{MJT20})} \label{cor: finiteness for homog spaces}
Assume that $G$ is linear, and let $X=G/H$ be a homogeneous space.
Then $X$ admits a finite number of equivalence classes of $(G,\sigma)$-equivariant real structures.
\end{corollary}

Indeed, if $\Aut^G(X)$ is a linear algebraic group, then it follows from \cite[Section 6.2, Th\'eor\`eme]{BS64} that $\H^1(\Gamma,\Aut^G(X))$ is finite, and this is in particular the case when $X \simeq G/H$ is homogeneous under the action of a linear algebraic group $G$ because then $\Aut^G(X) \simeq N_G(H)/H$ (see \cite[Proposition 1.8]{Tim11}).

\section{Real structures on linear representations and nilpotent orbits}\label{sec: warm-up}
In this section, we consider the Key Problem for two classical families of varieties with a reductive group action and related to representation theory, namely the linear representations (in Section \ref{subsec:real structures on linear reps}) and the nilpotent orbits in semisimple Lie algebras, as well as the normalizations of their closures (in Section \ref{subsec:real structures on nilpotent orbits}).

\subsection{Real structures on linear representations}\label{subsec:real structures on linear reps}
In this subsection, based on the appendix of \cite{MJT18} written by Borovoi (and based on the work of Tits \cite{Tits71}), we consider the Key Problem for finite-dimensional linear representations of complex simple groups. 
To simplify the situation, we only consider the case of irreducible representations.

\smallskip

Let $G$ be a complex (always assumed to be simply-connected) simple group, let $B \subseteq G$ be a Borel subgroup, let $U$ be the unipotent radical of $B$, and let $T \subseteq B$ be a maximal torus. Recall that the restriction map $\X(B) \to \X(T),\ \chi \to \chi_{|T}$ induces an isomorphism of character lattices, and so we can identify $\X=\X(T)$ with $\X(B)$. 
A \emph{dominant weight} of $G$ with respect to $(B,T)$ is a character $\lambda \in \X$ such that $\langle \lambda,\gamma \rangle \geq 0$ for every $\gamma \in \SS^\vee$. The subset $\D \subseteq \X$ of dominant weights satisfies $\D=\X \cap C$ in $\X_\Q$, where $C \subseteq \X_\Q$ is a polyhedral cone known as the \emph{fundamental Weyl chamber}.
\begin{example}
\begin{multicols}{2}
Let $G=\SL_3$ and let $\SS=\{\alpha_1,\alpha_2\}$ be the corresponding set of simple roots. The character lattice $\X$ is a rank two lattice generated by the fundamental weights $\varpi_1 = \frac23
  \alpha_1 + \frac13 \alpha_2$ and $\varpi_2 = \frac13 \alpha_1 + \frac23
  \alpha_2$, and the fundamental Weyl chamber is the grey cone. The dominant weights of $G$ are then all the nonnegative linear combinations of $\varpi_1$ and $\varpi_2$.
  \begin{center}
\includegraphics[width=4.5cm]{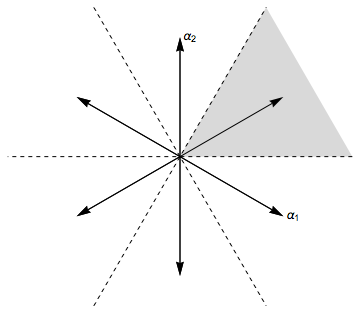}
\end{center}
\end{multicols}
\end{example}
Let $V$ be an irreducible representation of $G$. Then the fixed point set $V^U$ is a line on which $B$ acts through some character $\lambda \in \X$ called \emph{highest weight} of the representation $V$. The map that associates to an irreducible representation its highest weight induces a one-to-one correspondence between the irreducible representations (up to isomorphism) and the dominant weights. For a given $\lambda \in \D$, we denote by $V_\lambda$ the corresponding irreducible representation.

\smallskip

Let $\sigma$ be a real group structure on $G$. 
Let $Z=Z(G)$ be the center of $G$; it is stable by $\sigma$, and so we can view $Z$ as an abelian $\Gamma$-group and consider the abelian group $\H^2(\Gamma,Z)$. We denote by $\delta([\sigma]) \in  \H^2(\Gamma,Z)$ the \emph{Tits class} of $\sigma$ (see \cite[Section 4.2]{Tits71} or \cite[Section 1.3]{MJT18} for the definition of the Tits class); it is a cohomology class obtained from $\sigma$, whose conjugacy class by inner automorphisms of $G$ identifies with an element $[\sigma]$ of $\H^1(\Gamma,G/Z)$, 
 by applying the Galois cohomology functor to the short exact sequence of $\Gamma$-groups
\[
1 \to Z \to G \to G/Z \to 1.
\]
\begin{remark}
Tits classes for all complex simply-connected simple groups have been computed by Borovoi using Theorem \ref{th: real structures on linear reps} below; complete tables can be found in \cite[Appendix]{MJT18}.
\end{remark}

The $\star$-action induced by $\sigma$ stabilizes the subset $\D$ in $\X$  (see \cite[Lemma A.1]{MJT18}).
Let $\lambda \in \D$ such that $\gamma \cdot \lambda=\lambda$ (it will become clear with Theorem \ref{th: real structures on linear reps} why we only consider such dominant weights). This implies that the restriction $\lambda_{|Z}\colon Z \to \G_{m,\C}$ satisfies $\lambda_{|Z} \circ \sigma_{|Z}=\sigma_0 \circ \lambda_{|Z}$, with $\sigma_0(t)=\overline{t}$; in other words, $\lambda_{|Z}\colon Z \to \G_{m,\C}$ is a morphism of abelian  $\Gamma$-groups, and so it induces a homomorphism between second cohomology groups 
\[
\lambda_*\colon \H^2(\Gamma,Z) \to \H^2(\Gamma,\G_{m,\C})=\{\pm 1\}.
\]

If we restrict ourselves to considering only real forms in the category of linear representations, then the next result gives a complete answer to Question \ref{item: existence part of Key-Prob3} of the Key Problem for irreducible linear representations.

\begin{theorem} \emph{(\cite[Theorem A.2]{MJT18}, see also \cite[Theorem 7.2]{Tits71})}\\ \label{th: real structures on linear reps}\noindent We keep the previous notation, and we denote by $F=G/\langle \sigma \rangle$ the real form of $G$ corresponding to $\sigma$.
Let $V_\lambda$ be an irreducible representation of $G$.
Then $V_\lambda$ admits a real $F$-form (in the category of linear representations of $F$) if and only if 
\begin{enumerate}
\item \label{item i: th linear reps} the dominant weight $\lambda$ is fixed by the $\star$-action induced by $\sigma$; and
\item \label{item ii: th linear reps}  the cohomology class $\lambda_*(\delta([\sigma])) \in \H^2(\Gamma,\G_{m,\C})$ is trivial.
\end{enumerate}
\end{theorem}

\begin{remark}\item
\begin{itemize}
\item 
For any $\lambda \in \D$, the linear representation $V_\lambda \oplus V_{\gamma \cdot \lambda}$ admits a real $F$-form (in the category of linear representations of $F$). In particular, if only condition \ref{item i: th linear reps} of Theorem \ref{th: real structures on linear reps} holds, then $V_\lambda$ does not admit a real $F$-form, but $V_\lambda \oplus V_\lambda$ does (it is given by $V_\lambda$ itself, seen as a real vector space).
\item
If a linear representation of $G$ admits a real $F$-form (in the category of linear representations of $F$), then it follows from  Schur's lemma that this real $F$-form is unique (in the category of linear representations of $F$), up to isomorphism.
Tables with the list of all the irreducible representations of $G$ admitting a real $F$-form, in the category of linear representations of $F$, can be found in \cite{Tits67}.
\item
A linear representation of $G$ may admit a real $F$-form that is not a linear representation of $F$; see the last paragraph of Section \ref{subsec: real structures on T-varieties} for some references.
\end{itemize}
\end{remark}

\subsection{Real structures on nilpotent orbits in semisimple Lie algebras}\label{subsec:real structures on nilpotent orbits}
In this subsection, based on the work of Moser-Jauslin--Terpereau \cite{MJT21}, we consider the Key Problem for nilpotent orbits in complex semisimple Lie algebras and for the normalizations of their closures; these are examples of varieties with symplectic singularities (whose symplectic desingularizations are quite well understood), and they furthermore have a series of applications in representation theory of algebraic groups, Lie algebras and related objects (such as Weyl groups).
We refer to \cite{CM93} for a general reference on nilpotent orbits.

\smallskip

Let $G$ be a complex semisimple group, let $\sigma$ be a real group structure on $G$, and let $\g$ be the Lie algebra of $G$.
Recall that an element $x \in \g$ is called \emph{nilpotent} if the endomorphism $\mathrm{ad}_x\colon \g \to \g, \ y \mapsto [x,y]$ is nilpotent. Also, $G$ naturally acts on $\g$ via the adjoint action and nilpotent elements of $\g$ form a finite union of $G$-orbits called \emph{nilpotent orbits} in $\g$.
Denote by 
\[d\sigma_e\colon \g \to \g\] 
the differential of $\sigma\colon G \to G$ at the identity element (the fact that $\sigma$ fixes the identity element $e \in G$ follows directly from the definition of a real group structure). One can check that $d\sigma_e$ is a $(G,\sigma)$-equivariant real structure on $\g$, viewed as a $G$-variety for the adjoint action, and that $d\sigma_e$ stabilizes the (finite) set of nilpotent orbits in $\g$. 
In particular, if $\O$ is a nilpotent orbit in $\g$ such that $d\sigma_e(\O)=\O$, then $(d\sigma_e)_{|\O}$ is a $(G,\sigma)$-equivariant real structure on $\O$.
However, there are also equivariant real structures on nilpotent orbits that are not obtained by differentiating a real group structure on $G$ (see Example \ref{ex:non-gloabl1}), nor even by restricting an equivariant real structure from the Lie algebra $\g$ (see Example \ref{ex:non-gloabl2}).

\begin{example}\label{ex:non-gloabl1}
Let $\g$ be a semisimple Lie algebra, and assume that $d\sigma_e(\O)=\O$. Let $\theta \in \R$.
Then $\mu_\theta \colon \O \to \O,\ v \mapsto e^{i\theta} d\sigma_e(v)$ is a $(G,\sigma)$-equivariant real structure on $\O$ which cannot be obtained by differentiating a real group structure on $G$ when $\theta \notin 2\pi \Z$ (because, in this case, $\mu_\theta$ does not preserve the Lie bracket).
\end{example}

\begin{example}\label{ex:non-gloabl2}(see \cite[Example 1.4 and Section 6]{MJT21})
Let $G=\SL_3(\C)$ with $\sigma(g)=\overline{g}$ for all $g \in G$ (here $\overline{g}$ denotes the complex conjugate of $g$), and let $\O_{reg}$ be the \emph{regular} nilpotent orbit in $\sl_3$ (i.e.~the unique nilpotent orbit whose closure contains all the other nilpotent orbits).
A representative of $\O_{reg}$ is given by \small$\begin{bmatrix}
0 & 1 & 0 \\ 0 & 0 & 1\\ 0 &0&0
\end{bmatrix}$\normalsize; its stabilizer is
\small
$
H:=\Stab_G\left(\begin{bmatrix}
0 & 1 & 0\\ 0 & 0 & 1\\0 & 0 & 0
\end{bmatrix}\right)=\left\{ \begin{bmatrix}
a & b & c\\ 0 & a & b\\0 & 0 & a
\end{bmatrix};\ a,b,c \in \C \ \text{with}\ a^3=1\right\}$.
\normalsize
Then the map $\mu$ defined by \small 
\[
\O_{reg} \to \O_{reg},\ g \cdot \begin{bmatrix}
0 & 1 & 0 \\ 0 & 0 & 1\\ 0 &0&0
\end{bmatrix} \mapsto \left( \sigma(g)  \begin{bmatrix}
1 & i & 0 \\ 0 & 1 & 0\\ 0 &0&1
\end{bmatrix} \right) \cdot \begin{bmatrix}
0 & 1 & 0 \\ 0 & 0 & 1\\ 0 &0&0
\end{bmatrix} =\sigma(g) \cdot \begin{bmatrix}
0 & 1 & i \\ 0 & 0 & 1\\ 0 &0&0
\end{bmatrix}  
\] \normalsize
is a $(G,\sigma)$-equivariant real structure on $\O_{reg}$ that does not extend to a $(G,\sigma)$-equivariant real structure on $\sl_3$.
\end{example}

We recall that to every nilpotent orbit $\O \subseteq \g$ is associated a unique \emph{weighted Dynkin diagram} $w(\O)$, i.e.~a Dynkin diagram of type $\g$ where each node is labeled with $0$, $1$ or $2$; see \cite[Section 3.5]{CM93} for details. Let us note however that not every weighted Dynkin diagram corresponds to a nilpotent orbit of $\g$.
Also, denoting by $\gamma \cdot w(\O)$ the weighted Dynkin diagram obtained from $w(\O)$ by applying the $\star$-action induced by $\sigma$ (see Definition \ref{def: star-action}), \cite[Proposition 3.2]{MJT21} yields the equality
\[\gamma \cdot w(\O)=w(d\sigma_e(\O)).\] 
Alternatively, nilpotent orbits in complex simple Lie algebras can also be classified via partitions; see \cite[Section 5.1]{CM93} for details. Both classifications can be useful depending on the situation.
\begin{example}
Let $\g=\sl_3$. Then the nilpotent cone of $\g$ is the union of three nilpotent orbits, corresponding to the partitions $[3]$, $[2,1]$, and $[1^3]$ associated with the Jordan normal forms \small $\begin{bmatrix}
0 & 1& 0\\ 0 & 0 & 1\\ 0 & 0 & 0 
\end{bmatrix}$, $\begin{bmatrix}
0 & 1& 0\\ 0 & 0 & 0\\ 0 & 0 & 0 
\end{bmatrix}$, \normalsize and $\begin{bmatrix} 0
\end{bmatrix}$ respectively, and whose corresponding weighted Dynkin diagrams are
\begin{center}
  \begin{tikzpicture}[scale=.4]
    \draw (-1,0) node[anchor=east]{} ;
    \draw[xshift=0 cm,thick] (0 cm,0) circle (.2cm) node[anchor=south]{2};
      \draw[xshift=1 cm,thick] (1 cm,0) circle (.2cm) node[anchor=south]{2};
        \draw[xshift=2 cm,thick] (2 cm,0) circle (.2cm)node[anchor=south]{2};
    \draw[thick] (0.2 cm,0) -- +(1.6 cm,0);
     \draw[thick] (2.2 cm,0) -- +(1.6 cm,0);
     \end{tikzpicture}
       \begin{tikzpicture}[scale=.4]
    \draw (-1,0) node[anchor=east]{} ;
    \draw[xshift=0 cm,thick] (0 cm,0) circle (.2cm) node[anchor=south]{1};
      \draw[xshift=1 cm,thick] (1 cm,0) circle (.2cm) node[anchor=south]{0};
        \draw[xshift=2 cm,thick] (2 cm,0) circle (.2cm)node[anchor=south]{1};
    \draw[thick] (0.2 cm,0) -- +(1.6 cm,0);
     \draw[thick] (2.2 cm,0) -- +(1.6 cm,0);
     \end{tikzpicture}
       \begin{tikzpicture}[scale=.4]
    \draw (-1,0) node[anchor=east]{} ;
    \draw[xshift=0 cm,thick] (0 cm,0) circle (.2cm) node[anchor=south]{0};
      \draw[xshift=1 cm,thick] (1 cm,0) circle (.2cm) node[anchor=south]{0};
        \draw[xshift=2 cm,thick] (2 cm,0) circle (.2cm)node[anchor=south]{0};
    \draw[thick] (0.2 cm,0) -- +(1.6 cm,0);
     \draw[thick] (2.2 cm,0) -- +(1.6 cm,0);
     \end{tikzpicture}
\end{center}
\end{example}

The following theorem provides a complete answer to the Key Problem for nilpotent orbits and for the normalizations of their closures in $\g$.

\begin{theorem}\emph{(\cite[Main Theorem]{MJT21})}\label{th:main th nilpotent}
We keep the previous notation.
Let $\overline{\O}$ be the closure of $\O$ in $\g$, and let $\widetilde{\O}$ be the normalization of $\overline{\O}$.
Then the following are equivalent:
\begin{enumerate}[$(a)$]
\item\label{item:a} $\O$ is $d\sigma_e$-stable (i.e.~$d\sigma_e(\O)=\O$); 
\item\label{item:b} $\O$ admits a $(G,\sigma)$-equivariant real structure;
\item\label{item:c} $\widetilde{\O}$ admits a $(G,\sigma)$-equivariant real structure; and
\item\label{item:d} $\gamma \cdot w(\O)=w(\O)$.
\end{enumerate}
Moreover, if these equivalent assertions hold, then all $(G,\sigma)$-equivariant real structures on $\O$, resp.~on $\widetilde{\O}$, are equivalent. 
\end{theorem}

\begin{remark}\ 
\begin{itemize}

\item 
Even if nilpotent orbits may have equivariant real structures which cannot be obtained by differentiating a real group structure (see Example \ref{ex:non-gloabl1}), the important punchline of Theorem \ref{th:main th nilpotent} is that for answering Question \ref{item: existence part of Key-Prob3} of the Key Problem, it suffices to consider equivariant real structures that are obtained by differentiating real group structures.	

\item Nilpotent orbits that admit a $(G,\sigma)$-equivariant real structure are described in \cite[Section 3]{MJT21}.
It turns out that, except for a few cases in type $D_{2n}$ (with $n \geq 2$), every nilpotent orbit $\O$ in $\g$  admits a $(G,\sigma)$-equivariant real structure when $\g$ is simple (see Example \ref{ex:D4} below for the $D_4$ case).

\item  When $\overline{\O}$ is non-normal, we do not know whether every $(G,\sigma)$-equivariant real structure on $\O$ extends to $\overline{\O}$.
A brief review on what is known about the \mbox{(non-)}normality of nilpotent orbit closures in semisimple Lie algebras can be found in \cite[Section 5]{MJT21}.  

\item Assume that $G$ is simple, and let $F=G/\langle \sigma \rangle$ be the real form corresponding to $\sigma$. For every nilpotent orbit $\O$ in $\g$, we have that $\O(\C)^{d\sigma_e}=\O(\C) \cap \g^{d\sigma_e}$ is a real manifold (possibly empty) whose $F(\R)$-orbits, usually called \emph{real nilpotent orbits}, are classified (see \cite[Section 9]{CM93}). 
\end{itemize}
\end{remark}

\begin{example}(\cite[Example 3.6]{MJT21})\label{ex:D4}
A partial order relation is defined on the set of nilpotent orbits of $\g$ by inclusion of closures.
In type $D_4$, the Hasse diagram representing this partial order is the following (see \cite[Section 6.2]{CM93}).
\begin{multicols}{2}
\begin{center}
\scalebox{0.8}{
\xymatrix@R=4mm@C=2cm{
& \O_{[7,1]} \ar@{-}[d] & \\
& \O_{[5,3]} \ar@{-}[rd] \ar@{-}[ld] \ar@{-}[d] &\\
\O_{[4^2]}^{\mathrm{I}} \ar@{-}[rd] \ar@/^3pc/@{<.>}[rr] \ar@{<.>}[r] & \O_{[5,1^3]} \ar@{<.>}[r]\ar@{-}[d]& \O_{[4^2]}^{\mathrm{I\!I}} \ar@{-}[ld]\\
& \O_{[3^2,1^2]} \ar@{-}[d]&  \\
& \O_{[3,2^2,1]} \ar@{-}[ld] \ar@{-}[rd] \ar@{-}[d]&  \\
\O_{[2^4]}^{\mathrm{I}} \ar@/^3pc/@{<.>}[rr] \ar@{<.>}[r]\ar@{-}[rd] & \O_{[3,1^5]} \ar@{<.>}[r]\ar@{-}[d]& \O_{[2^4]}^{\mathrm{I\!I}}  \ar@{-}[ld] \\  
& \O_{[2^2,1^4]}\ar@{-}[d]&  \\
& \O_{[1^8]}&   
      }
            }
\end{center}
Here the dotted arrows indicate which pairs of orbits can be swapped by the $\star$-action induced by an outer twist of a split real group structure on $G$.
\end{multicols}
\end{example}

\section{Real structures on varieties with a torus action} \label{sec: real structures on T-varieties}
In this section, we consider the Key Problem when the (effectively acting) reductive group is a complex torus $T$. There are basically two cases to consider: 
\begin{itemize}
\item The first case is when $T$ acts with a dense open orbit (case of \emph{toric varieties}), in which case there is a well-known combinatorial description of the $T$-varieties in terms of fans, and a complete answer to the Key Problem can be provided. 
\item The second case is when a general $T$-orbit is of codimension $c\geq 1$ (this number $c$ is called the \emph{complexity} of the $T$-variety), in which case there is again a combinatorial description due to Altmann--Hausen \cite{AH06} (affine case) and Altmann--Hausen--S\"u\ss \ \cite{AHS08} (general case) that makes it possible to give a satisfactory answer to Question \ref{item: existence part of Key-Prob3} of the Key Problem, at least in the affine case, but at the cost of a rather heavy technical machinery. (Let us mention that the combinatorial description of Altmann-Hausen extends combinatorial descriptions obtained by Demazure in \cite{Dem88} and by Flenner--Zaidenberg in \cite{FZ03}.)
\end{itemize}

\subsection{Real structures on toric varieties}\label{subsec: real structures on toric varieties}
We refer to \cite{Ful93} or \cite{CLS11} for a complete account on toric varieties over the field of complex numbers.
Let $T \simeq \G_{m,\C}^n$ be a complex torus of dimension $n \geq 1$, let $\sigma$ be a real group structure on $T$, and let $X$ be a complex  toric variety with open orbit $X_0 \iso T$. Our goal in this subsection is to explain how to determine all the equivalence classes of $(T,\sigma)$-equivariant real structures on $X$. 

Let us note that the Galois descent $\C/\R$ is always effective for complex toric varieties. Indeed, by Sumihiro's theorem (\cite[Corollary 2]{Sum74}), a toric variety is covered by $T$-stable affine open subsets, which are permuted by the $\Gamma$-action, and the union of two $T$-stable affine open subsets is quasi-projective by \cite[Corollary 3.2.12]{Per14} (see also \cite[Proposition 1.9 and Remark 1.11]{Hur11}). Hence, there is a one-to-one correspondence between equivalence classes of equivariant real structures on $X$ and isomorphism classes of real forms of $X$ in the category of real toric varieties. 

\smallskip

By Lemma \ref{lem: real form on tori}, there exists $\psi \in \Aut_{\gr}(T)$ and $n_0,n_1,n_2 \in \N_{\geq 0}$ such that 
\[\sigma=\psi \circ (\sigma_0^{\times n_0}\times\sigma_1^{\times n_1}\times\sigma_2^{\times n_2}) \circ \psi^{-1},\ \ \text{with}\ n_0+n_1+2n_2=n.\] 
Furthermore, denoting 
\[\tau_1: \G_{m,\C} \to \G_{m,\C},\ t \mapsto -\overline{{t}^{-1}},\] 
which is a $(\G_{m,\C},\sigma_1)$-equivariant real structure on $\G_{m,\C}$, we check that the antiregular involution $X_0 \to X_0$ defined by
\[
\psi \circ (\sigma_0^{\times n_0}\times \rho_1 \times \cdots \times \rho_{n_1} \times\sigma_2^{\times n_2}) \circ \psi^{-1}, \ \text{where}\ \rho_i\in \{ \sigma_1, \tau_1\} \text{ for each } i=1,\ldots, n_1, 
\]
is a $(T,\sigma)$-equivariant real structure on $X_0 \simeq T$ (seen as a $T$-variety for the usual multiplication on $T$). Furthermore, any $(T,\sigma)$-equivariant real structure on $X_0$ is equivalent to one of this kind, and two equivariant real structures $\psi \circ (\sigma_0^{\times n_0}\times \rho_1 \times \cdots \times \rho_{n_1} \times\sigma_2^{\times n_2}) \circ \psi^{-1}$ and $\psi \circ (\sigma_0^{\times n_0}\times \rho'_1 \times \cdots \times \rho'_{n_1} \times\sigma_2^{\times n_2}) \circ \psi^{-1}$ on $X_0$ are equivalent if and only if $\rho_i=\rho'_i$ for each $i=1,\ldots, n_1$; in particular, $X_0$ admits exactly $2^{n_1}$ equivalence classes of $(T,\sigma)$-equivariant real structures. (This can be proved for instance using Proposition \ref{prop:Galois H1 to param eq real structures} since $\H^1(\Gamma,\Aut^{T}(X_0)) \simeq (\Z/2\Z)^{n_1}$ by \cite[Proposition 1.18]{MJT18}.)

\begin{remark}
Let $\mu$ be a $(T,\sigma)$-equivariant real structure on $X_0$. Then, either $\mu$ is equivalent to $\sigma$, in which case $X_0/\langle \mu \rangle \simeq T/\langle \sigma \rangle$ is  the trivial $(T/\langle \sigma \rangle)$-torsor, or $X_0/\langle \mu \rangle$ is a $(T/\langle \sigma \rangle)$-torsor with no real points (in particular it is irrational).
\end{remark}

\smallskip

It remains to determine which $(T,\sigma)$-equivariant real structures on $X_0$ extend to $X$, and which ones stay equivalent after extension. Let us first note that $\Aut^T(X)\simeq \Aut^T(X_0) \simeq T$, which implies two things:
\begin{itemize}
\item If $\mu_1$ and $\mu_2$ are two equivalent $(T,\sigma)$-equivariant real structures on $X_0$, then $\mu_1$ extends to $X$ if and only if $\mu_2$ extends to $X$.
\item Two $(T,\sigma)$-equivariant real structures on $X$ are equivalent if and only if their restrictions on $X_0$ are equivalent. 
\end{itemize}
Therefore, we are left with determining which $(T,\sigma)$-equivariant real structures $\psi \circ (\sigma_0^{\times n_0}\times \rho_1 \times \cdots \times \rho_{n_1} \times\sigma_2^{\times n_2}) \circ \psi^{-1}$ on $X_0$ extend to $X$. It turns out that the answer, obtained by Huruguen in \cite{Hur11},  only depends on $\sigma$.

\begin{theorem}\label{th:real structure extend to toric embeddings} \emph{(\cite[Proposition 1.19 and Theorem 1.25]{Hur11})}\\
With the previous notation, a $(T,\sigma)$-equivariant real structure $\mu$ on $X_0$ extends to $X$ if and only if the $\star$-action on $\X_\Q^\vee$ induced by $\sigma$  (see Example \ref{ex: star-action for tori}) stabilizes the fan $\Sigma$ associated with the $T$-equivariant embedding $X_0 \hookrightarrow X$.
\end{theorem}

It follows from Theorem \ref{th:real structure extend to toric embeddings} and the discussion above that the complex toric variety $X$ admits either $2^{n_1}$ equivalence classes of $(T,\sigma)$-equivariant real structures (when the $\star$-action on $\X_\Q^\vee$ stabilizes the fan $\Sigma$) or none. In particular, we have a complete answer to the Key Problem for toric varieties in terms of combinatorial objects easy to implement in a software system. We now give some examples.

\begin{example}
Let $T \simeq \G_{m,\C}$, and let $\Sigma$ be the fan in $\X_\Q^\vee \simeq \Q$ generated by $e_1:=1$; it corresponds to the natural embedding $X_0=\G_{m,\C} \hookrightarrow X=\A_\C^1$. 
\begin{itemize}
\item If $\sigma=\sigma_0$, then the corresponding $\star$-action on $\X_\Q^\vee$ is trivial, and so any $(T,\sigma_0)$-equivariant real structure on $X_0$ extends to $\A_\C^1$ (and they are all equivalent); it corresponds to the real $\G_{m,\R}$-form $\A_\R^1$, up to isomorphism. 
\item If $\sigma=\sigma_1$, then the corresponding $\star$-action on $\X_\Q^\vee$ is $r \mapsto -r$, and so neither of the two inequivalent $(T,\sigma_1)$-equivariant real structures $\sigma_1$ and $\tau_1$ extends to $\A_\C^1$, i.e.~$\A_\C^1$ has no real $\Sone$-form. 
\end{itemize}
However, if we now take $\Sigma$ the be the complete fan in $\X_\Q^\vee \simeq \Q$ generated by $e_1$ and $-e_1$, then any $(T,\sigma)$-equivariant real structure on $X_0$ extends to $X=\P_\C^1$ since both $\star$-actions stabilize $\Sigma$; it gives rise to the real $\G_{m,\R}$-form $\P_{\R}^{1}$ when $\sigma=\sigma_0$, and to the two real $\Sone$-forms $\P_\R^1$ and the empty conic when $\sigma=\sigma_1$.
\end{example}

Let us mention that if $\sigma$ and $\sigma'$ are two equivalent real group structures on $T$, it may happen that $X$ admits $2^{n_1}$ inequivalent $(T,\sigma)$-equivariant real structures but no $(T,\sigma')$-equivariant real structure. Hence, it is essential to set a real group structure $\sigma$, and not just its equivalence class, when studying the equivariant real structures on toric varieties.

\begin{example}
Let $T \simeq \G_{m,\C}^2$, and let $\sigma$ and $\sigma'$ be the two equivalent real group structures defined in Example \ref{ex:counter-ex not strongly eq}. The $\star$-action on $\X^\vee \simeq \Z^2$ induced by $\sigma$ is given by $(m,n) \mapsto (m,-n)$, while the one induced by $\sigma'$ is given by $(m,n) \mapsto (m-2n,-n)$. Consider now the $T$-equivariant embedding $X_0 \hookrightarrow X \simeq \P_\C^1 \times \P_\C^1$ corresponding to the following complete fan in $\X_\Q^\vee$: 
\begin{center}
\begin{tikzpicture} 
\draw[->] (5,0) -- (6,0); 
\draw[->] (5,0) -- (4,0); 
\draw[->] (5,0) -- (5,1); 
\draw[->] (5,0) -- (5,-1); 
\node at (6.3,0) {$e_1$}; 
\node at (5,1.2) {$e_2$}; 
\node at (3.6,0) {$-e_1$}; 
\node at (5,-1.2) {$-e_2$}; 
\node at (5.6,0.6) {$C_1$}; 
\node at (4.4,0.6) {$C_2$}; 
\node at (5.6,-0.6) {$C_4$}; 
\node at (4.4,-0.6) {$C_3$}; 
\end{tikzpicture}
\end{center}

\noindent We denote this fan by $\Sigma$. 
The $\star$-action induced by $\sigma$ swaps the cones $C_1 \leftrightarrow C_4$ and $C_2 \leftrightarrow C_3$, so $\Sigma$ is $\Gamma$-stable and $X$ admits $2$ equivalence classes of $(T,\sigma)$-equivariant real structures (corresponding to the isomorphism classes of the two real $(\G_{m,\R} \times \Sone)$-forms $\P_\R^1 \times \P_\R^1$ and $\varnothing$). But the $\star$-action induced by $\sigma'$ maps $e_2$ to $-2e_1-e_2$, so $\Sigma$ is not $\Gamma$-stable and $X$ admits no $(T,\sigma')$-equivariant real structure.
\end{example}

\begin{works}\ 
\begin{itemize}
\smallskip
\item Equivariant real structures (in a slightly broader sense than in the present article)
on projective smooth toric surfaces and threefolds have been classified by Delaunay in \cite{Del03,Del04}; she has also studied the topology of the corresponding real loci.
\smallskip
\item Descent data and $\K$-forms for toric varieties, with $\K$ an arbitrary field, have been studied by Huruguen in \cite{Hur11} for arbitrary toric varieties (focusing on determining when a given descent datum on the open orbit extends to the whole variety), by Elizondo--Lima-Filho--Sottile--Teitler in \cite{ELST14} for projective spaces and toric surfaces (with a detailed classification of the possible descent data in these cases), and by Duncan in \cite{Dun16} who has considered and compared $\K$-forms of toric varieties in three different categories, namely the abstract varieties, the toric varieties, and the \emph{neutral} toric varieties (i.e.~toric varieties whose open orbit is a trivial torsor).
\end{itemize}
\end{works}

\subsection{Real structures on arbitrary affine \texorpdfstring{$T$}{T}-varieties}\label{subsec: real structures on T-varieties}
In this subsection, based on the work of Gillard \cite{Gil20,Gil23}, we still consider torus actions on complex algebraic varieties, but we do not assume that there is an open orbit anymore, and we restrict ourselves to the affine setting (in which case the Galois descent $\C/\R$ is always effective). Our main goal in this subsection is to outline the combinatorial answer obtained by Gillard to Question \ref{item: existence part of Key-Prob3} of the Key Problem for affine $T$-varieties.

\smallskip

We start by recalling, very briefly, the combinatorial description of affine $T$-varieties due to Altmann--Hausen (see \cite{AH06} for details):
Let $T \simeq \G_{m,\C}^n$ be a complex torus, with character lattice $\X \simeq \Z^n$, and let $X$ be a complex affine $T$-variety. The $T$-action on $X$ corresponds to an $\X$-grading of the coordinate ring $\C[X]$ as follows:
\[\C[X] = \bigoplus_{\chi \in \X} \C[X]_\chi \text{ with } \C[X]_\chi:=\{ f \in \C[X] \ |\ \forall t \in T, \forall x \in X,\  (t^{-1} \cdot f)(x)=\chi(t)f(x)\}.\] 
Let $\omega$ be a full dimensional polyhedral cone in $\X_\Q$, let $Y$ be a \textsl{semiprojective} variety (i.e.~$\H^0(Y,\O_Y)$ is a finitely generated $\C$-algebra and the affinization morphism $Y \to\Spec(H^0(Y,\O_Y))$ is projective), and let $\DDD:= \sum \Delta_i \otimes D_i$ be a  \textsl{polyhedral divisor} on $Y$; this means  that the $D_i$ are prime divisors on $Y$ and the coefficients $\Delta_i$ are convex polyhedra in $\X_\Q^\vee$ that can be written as the Minkowski sum of a polytope and the cone $\omega$ (we then say that $\omega$ is the \emph{tail cone} of $\Delta_i$).
Then, for every $\chi \in \omega \cap \X$, a Weil $\Q$-divisor on $Y$ is given by
$\DDD(\chi):= \sum \text{min}\{ \langle \chi | \Delta_i \rangle \} \otimes D_i$, where $\langle -|- \rangle\colon \X_\Q \times \X_\Q^\vee \to \Q$ denotes the dual pairing. 
Assume furthermore that $\DDD$ is \emph{proper}, i.e. that any $\DDD(\chi)$ is a semiample rational Cartier divisor on $Y$, which is big whenever $\chi$ belongs to the relative interior of the cone $\omega$.
From such a triple ($\omega, Y, \DDD$), that we will call an \emph{AH datum} in the rest of this subsection, Altmann and Hausen construct an $\X$-graded $\C$-algebra defined by 
\[
A[\omega, Y, \DDD]:= \bigoplus_{\chi \in \omega \cap \X} \H^0(Y, \mathcal{O}_Y(\DDD(\chi))) \subseteq \C(Y)[\X].
\]  
\begin{theorem} \emph{(see \cite[Theorems 3.1 and 3.4]{AH06})}\label{th:AH complex affine}
\begin{enumerate}
\item Let $(\omega, Y, \DDD)$ be an AH datum. Then $Z[\omega,Y, \DDD] := \Spec(A[\omega,Y, \DDD])$  is a complex affine $T$-variety that is $T$-equivariantly birational to $T \times Y$.
\item\label{item:ii of th 3.4} Let $X$ be a complex affine $T$-variety. Then there exists an AH datum  $(\omega, Y, \DDD)$ such that there is an isomorphism of $T$-varieties $X \simeq Z[\omega, Y, \DDD]$.
\end{enumerate}
\end{theorem}

\begin{remark} 
The semiprojective variety $Y$ in Theorem \ref{th:AH complex affine} \ref{item:ii of th 3.4} is not unique for a given $X$; a model can be obtained as the normalization of the main component of the inverse limit over all GIT-quotients of $X$ or by \emph{toric downgrading} (see \cite[Section 11]{AH06} for details).
\end{remark}

Based on the Altmann--Hausen combinatorial description of affine $T$-varieties (Theorem \ref{th:AH complex affine}), Gillard has obtained the following criterion that provides a complete answer to Question \ref{item: existence part of Key-Prob3} of the Key Problem for affine $T$-varieties.

\begin{theorem}\emph{(\cite[Theorems~4.3 and~4.6]{Gil20})}\label{th:AH real affine}
Let $\sigma$ be a real group structure on $T$. 
\begin{enumerate}
\item\label{item: i of th 3.6} Let $(\omega, Y, \DDD)$ be an AH datum, and let $X:=Z[\omega,Y, \DDD]$ be the corresponding complex affine $T$-variety.
If there exists a real structure $\rho_Y$ on $Y$ and a monoid homomorphism $h\colon \omega \cap \X \to \C(Y)^*$ such that
\begin{equation}\label{eq:monoid homomorphism}
\forall \chi \in \omega \cap \X,\ \rho_Y^*(\DDD(\chi))=\DDD(\ga \chi)+\div_Y(h(\ga \chi)) \ \text{and}\ h(\chi)\rho_Y^*(h(\ga \chi))=1,
\end{equation} 
where $\Gamma$ acts on $\X$ through the $\star$-action induced by $\sigma$ (see Example \ref{ex: star-action for tori}), then $X$ admits a $(T,\sigma)$-equivariant real structure $\mu_{[\omega,(Y,\rho_{Y}),\DDD,h]}$.
\item\label{item: ii of th 3.6} Let $X$ be an affine $T$-variety that admits a $(T,\sigma)$-equivariant real structure. 
Then there exists an AH datum  $(\omega, Y, \DDD)$, together with a real structure $\rho_Y$ on $Y$ and a monoid homomorphism $h\colon \omega \cap \X \to \C(Y)^*$ satisfying \eqref{eq:monoid homomorphism}, such that there is a $\Gamma$-equivariant isomorphism of $T$-varieties $X \simeq Z[\omega,Y,\DDD]$. 
\end{enumerate}
\end{theorem}

\begin{remark}\ 
\begin{itemize}
\item The important punchline in Theorem \ref{th:AH real affine} is that, as a consequence of the non-uniqueness of the combinatorial description of affine $T$-varieties due to Altmann–Hausen (see \cite[Section 8]{AH06} or \cite[Section 3.2.5]{Gil23}), one cannot just pick any AH datum to make the Galois descent work, but there always exists an AH datum that does the job.
\item 
The monoid homomorphism $h$ that appears in Theorem \ref{th:AH real affine} is a cocycle  that encodes a certain torsor over the generic point of the real variety associated to $(Y,\rho_Y)$; see \cite[Remark 3.2.2]{Gil23} for details.
In particular, when the real torus associated to $(T,\sigma)$ has no nontrivial torsor (which is the case if and only if it contains no $\mathbb{S}^1$ factor), we can take $h=1$, and then condition \eqref{eq:monoid homomorphism} becomes
\begin{equation}
\forall \chi \in \omega \cap \X,\ \rho_Y^*(\DDD(\chi))=\DDD(\ga \chi).
\end{equation} 
\item The criterion obtained by Gillard is in fact constructive in the sense that it produces an explicit $(T,\sigma)$-equivariant real structure $\mu_{[\omega,(Y,\rho_{Y}),\DDD,h]}$ on $Z[\omega,Y,\DDD]$. 
\item For a given complex affine $T$-variety $X$ admitting a $(T,\sigma)$-equivariant real structure, Gillard gives an explicit recipe to construct the pair $(Y,\rho_Y)$ by \emph{$\Gamma$-equivariant toric downgrading} (see \cite[Section 4.1]{Gil20}). 
\end{itemize}
\end{remark}

Theorem \ref{th:AH real affine} \ref{item: ii of th 3.6} simplifies a bit when the real group structure $\sigma$ on $T$ is equivalent to  $\sigma_0^{\times n_0}\times\sigma_2^{\times n_2}$ (no $\sigma_1$ factor, with the notation of Lemma \ref{lem: real form on tori}); indeed, we can then find $\DDD$ on $Y$ such that $\rho_Y^*(\DDD(\chi))=\DDD(\ga \chi)$, i.e.~take $h=1$ (see \cite[Section 5]{Lan15}). However, this simplification is not always possible when there is a $\sigma_1$ factor (see Example \ref{ex 4.10} below).

\begin{example}{({see \cite[Examples 2.14, 3.6 and 5.3]{Gil20} for details, more examples can be found in \cite[Section 4.1]{Gil23}})}
Let $T \simeq \G_{m,\C}^2$ and let $\sigma:=\sigma_2$ (with the notation of Lemma \ref{lem: real form on tori}). Let $X:=\A_\C^4$ on which $T$ acts by
\[
(s,t)\cdot (x_1,x_2,x_3,x_4):=(sx_1,tx_2,st^2 x_3,s^2tx_4), \ 
\]
and let $\mu$ be the $(T,\sigma)$-equivariant real structure on $X$ defined by 
\[
\mu\colon\ X \to X,\ (x_1,x_2,x_3,x_4)\mapsto (\overline{x_2},\overline{x_1},\overline{x_4},\overline{x_3});
\]
it corresponds to a real $R_{\C/\R}(\G_{m,\C})$-variety isomorphic to $\A_\R^4$.
Let $\{e_1,e_2\}$ be the canonical basis of $\X \simeq \Z^2$ and let $\omega:=\Cone(e_1,e_2)$. Let $Y$ be the projective toric surface corresponding to the fan on the left below, and let $\DDD:=\sum_{i=1}^{4} \Delta_i \otimes D_i$ be the polyhedral divisor on $Y$ such that each divisor $D_i$ corresponds to the ray generated by $v_i$ and the coefficients $\Delta_i \subseteq \X_\Q^\vee$ are the convex polyhedra pictured on the right below.
\begin{multicols}{4}
\definecolor{ffqqqq}{rgb}{1.,0.,0.}
\scalebox{0.65}{
\begin{tikzpicture}[line cap=round,line join=round,>=triangle 45,x=1.0cm,y=1.0cm]
\begin{axis}[
x=1.0cm,y=1.0cm,
axis lines=middle,
ymajorgrids=true,
xmajorgrids=true,
xmin=-2.5,
xmax=1.5,
ymin=-2.25,
ymax=1.5,
xtick={-3.0,-2.0,...,3.0},
ytick={-3.0,-2.0,...,3.0},]
\clip(-3.,-3.) rectangle (3.,3.);
\draw [->,line width=1.pt] (0.,0.) -- (1.,0.);
\draw [->,line width=1.pt] (0.,0.) -- (0.,1.);
\draw [->,line width=1.pt] (0.,0.) -- (-2.,-1.);
\draw [->,line width=1.pt] (0.,0.) -- (-1.,-2.);
\begin{scriptsize}
\draw[color=black] (1.2,0.15) node {$v_1$};
\draw[color=black] (0.2,1.2) node {$v_2$};
\draw[color=black] (-2.25,-1) node {$v_3$};
\draw[color=black] (-1.5,-2) node {$v_4$};
\end{scriptsize}
\end{axis}
\end{tikzpicture}
}

\scalebox{0.7}{

\begin{tikzpicture}[line cap=round,line join=round,>=triangle 45,x=1.0cm,y=1.0cm]
\begin{axis}[
x=1.0cm,y=1.0cm,
axis lines=middle,
ymajorgrids=true,
xmajorgrids=true,
xmin=-0.2,
xmax=2.75,
ymin=-0.2,
ymax=2.75,
xtick={-3.0,-2.0,...,3.0},
ytick={-3.0,-2.0,...,3.0},]
\clip(-3.,-3.) rectangle (3.,3.);
\shade[top color = white, bottom color = red] (0.,0.) -- (0.,3.) -- (3.,3.) -- (3.,0.)   -- cycle;
\draw [line width=2.pt,color=ffqqqq] (0.,0.)-- (0.,3.);
\draw [line width=2.pt,color=ffqqqq] (0.,3.)-- (3.,3.);
\draw [line width=2.pt,color=ffqqqq] (3.,3.)-- (3.,0.);
\draw [line width=2.pt,color=ffqqqq] (3.,0.)-- (0.,0.);
\node[label=right:{$ {\Delta_1=\Delta_2 } $}] (n1) at (0.3,1.) {};
\end{axis}
\end{tikzpicture}
}

\scalebox{0.7}{
\begin{tikzpicture}[line cap=round,line join=round,>=triangle 45,x=1.0cm,y=1.0cm]
\begin{axis}[
x=1.0cm,y=1.0cm,
axis lines=middle,
ymajorgrids=true,
xmajorgrids=true,
xmin=-0.2,
xmax=2.75,
ymin=-0.2,
ymax=2.75,
xtick={-3.0,-2.0,...,3.0},
ytick={-3.0,-2.0,...,3.0},]
\clip(-3.,-3.) rectangle (3.,3.);
\shade[top color = white, bottom color = red](0.,1.) -- (0.,3.) -- (3.,3.) -- (3.,0.) -- (2.,0.)  -- cycle;
\draw [line width=2.pt,color=ffqqqq] (0.,1.)-- (0.,3.);
\draw [line width=2.pt,color=ffqqqq] (0.,3.)-- (3.,3.);
\draw [line width=2.pt,color=ffqqqq] (3.,3.)-- (3.,0.);
\draw [line width=2.pt,color=ffqqqq] (3.,0.)-- (2.,0.);
\draw [line width=2.pt,color=ffqqqq] (2.,0.)-- (0.,1.);
\node[label=right:{$ {\Delta_3 } $}] (n1) at (0.75,1.) {};
\end{axis}
\end{tikzpicture}
}

\definecolor{ffqqqq}{rgb}{1.,0.,0.}
\scalebox{0.7}{
\begin{tikzpicture}[line cap=round,line join=round,>=triangle 45,x=1.0cm,y=1.0cm]
\begin{axis}[
x=1.0cm,y=1.0cm,
axis lines=middle,
ymajorgrids=true,
xmajorgrids=true,
xmin=-0.2,
xmax=2.75,
ymin=-0.2,
ymax=2.75,
xtick={-3.0,-2.0,...,3.0},
ytick={-3.0,-2.0,...,3.0},]
\clip(-3.,-3.) rectangle (3.,3.);
\shade[top color = white, bottom color = red] (0.,2.) -- (0.,3.) -- (3.,3.) -- (3.,0.) -- (1.,0.)  -- cycle;
\draw [line width=2.pt,color=ffqqqq] (0.,2.)-- (0.,3.);
\draw [line width=2.pt,color=ffqqqq] (0.,3.)-- (3.,3.);
\draw [line width=2.pt,color=ffqqqq] (3.,3.)-- (3.,0.);
\draw [line width=2.pt,color=ffqqqq] (3.,0.)-- (1.,0.);
\draw [line width=2.pt,color=ffqqqq] (1.,0.)-- (0.,2.);
\node[label=right:{$ {\Delta_4 } $}] (n1) at (0.75,1.) {};
\end{axis}
\end{tikzpicture}
}

\end{multicols}
\noindent Then ($\omega,Y,\DDD$) is an AH datum such that $X \simeq Z[\omega,Y,\DDD]$. Moreover, if we denote by $\rho_Y$ the extension to $Y$ of the real group structure $\sigma_2$ on the dense open orbit $Y_0 \simeq \G_{m,\C}^2$, and take $h=1$, then one can check that the pair $(\rho_Y,h)$ satisfies \eqref{eq:monoid homomorphism}.
\end{example}

\begin{example}{(see \cite[Section 4.1]{DL22})} \label{ex 4.10}
Let $T \simeq \G_{m,\C}$, let $\sigma:=\sigma_1$ (the real group structure corresponding to the circle $\Sone$), and let $\omega=\X_\Q \simeq \Q$. 
Let $Y$ be the complex surface with the real structure $\rho_Y$ whose corresponding real form is the real sphere $\mathbb{S}^2:=\Spec(\R[x,y,z]/(x^2+y^2+z^2-1))$. 
Let $\DDD:=\{1\} \otimes D$, where $D:=\{ 1-z=x+iy=0\}$ is a Cartier divisor on $Y=\mathbb{S}_\C^2$, and let $h\colon \Z \to \C(Y)^*,\ n \mapsto (1-z)^{-n}$. 
Then $(\omega,Y,\DDD)$ is an AH datum, the pair $(\rho_Y,h)$ satisfies \eqref{eq:monoid homomorphism}, and the corresponding $(T,\sigma)$-equivariant real structure $\mu_{[\omega,(Y,\rho_{Y}),\DDD,h]}$ on $X:=Z[\omega,Y,\DDD]$ corresponds to the real $\mathbb{S}^1$-variety $\mathbb{S}^3$ whose quotient $f\colon \mathbb{S}^3 \to \mathbb{S}^3/\mathbb{S}^1 \simeq \mathbb{S}^2$ is an \emph{algebraic model} of the famous Hopf fibration (which means that the continuous map $\mathbb{S}^3(\R) \to  \mathbb{S}^2(\R)$ induced by $f$ is the Hopf fibration). Let us note that, in this example, it is not possible to reduce to the case $h=1$.
\end{example}

As mentioned already, Theorem \ref{th:AH real affine} provides a complete answer to Question \ref{item: existence part of Key-Prob3} of the Key Problem for affine $T$-varieties. It is then natural to ask what about an answer for Question \ref{item: quantity part of Key-Prob3} of the Key Problem for these varieties. Unfortunately, for an arbitrary affine $T$-variety $X$, the group $\Aut^T(X)$ can be quite wild, and so in general it is hard enough even to determine if the set $\H^1(\Gamma,\Aut^T(X))$, which parametrizes the equivalence classes of equivariant real structures on $X$, is finite or not.

\begin{works}\ 
\begin{itemize}
\smallskip
\item Descent data and $\K$-forms for affine complexity-one $T$-varieties, with $\K$ an arbitrary field, have been studied by Langlois in \cite{Lan15}; in particular, if $\K=\R$  and the real group structure $\sigma$ on $T$ is equivalent to  $\sigma_0^{\times n_0}\times\sigma_2^{\times n_2}$ (no $\sigma_1$ factor), then \cite[Theorem 5.10]{Lan15} coincides with Theorem \ref{th:AH real affine}.

\smallskip

\item Real $\Sone$-forms of $\G_{m,\C}$-varieties have been described by Dubouloz--Liendo in \cite{DL22}; in particular, \cite[Theorem 2.7]{DL22} coincides with Theorem \ref{th:AH real affine} when $T=\G_{m,\C}$ and $\sigma=\sigma_1$.
Also, the case of smooth $\G_{m,\C}$-surfaces has been considered more closely by Dubouloz--Petitjean in \cite{DP20}; as an application, they establish that every compact differentiable $\Sone(\R)$-surface admits a unique smooth rational real affine model up to $\Sone$-equivariant birational diffeomorphism \cite[Theorem 1]{DP20}.

\smallskip

\item An example of a non-linearizable $\Sone$-action on $\A_\R^4$, whose complexification is a linearizable $\G_{m,\C}$-action on $\A_\C^4$, was obtained by Freudenburg--Moser-Jauslin in \cite{FMJ04}. (Note that Koras--Russell proved in \cite{KR13} that any $\Sone$-action on $\A_\R^3$ is linearizable.) Later, a systematic approach to find uncountably many non-linearizable $\Sone$-actions on $\A_\R^4$ that are all pairwise inequivalent was obtained by Moser-Jauslin in \cite{MJ19} (see also \cite[Section 4.2]{DL22}).
\end{itemize}
\end{works}

\section{Real structures on spherical varieties}\label{sec: real structures on spherical varieties}
In this section, we leave the world of torus actions and consider the Key Problem for another famous class of algebraic varieties endowed with an algebraic group action, namely the \emph{spherical varieties}.

Let $G$ be a complex reductive group (of simply-connected type). A $G$-variety $X$ is called \emph{spherical} if a Borel subgroup of $G$ acts on $X$ with a dense open orbit. This condition implies of course that $G$ itself acts on $X$ with a dense open orbit, and consequently the combinatorial description of spherical varieties is generally done in two steps.
\begin{itemize}
\item First step: classify the spherical homogeneous spaces under the $G$-action. 
\item Second step: for a given spherical homogeneous space $G/H$, classify the $G$-equivariant embeddings of $G/H$ (also known as \emph{spherical embeddings}). 
\end{itemize}
We will not go into details regarding the combinatorics that appear at each of these two steps to avoid making this survey article too technical. We refer to \cite{Tim11} and the references therein for a complete account on spherical homogeneous spaces and their equivariant embeddings.

\subsection{Real structures on spherical homogeneous spaces}\label{subsec:real structures on spherical homog spaces}
In this subsection, we consider the Key Problem for spherical homogeneous spaces.
Let $G$ be a complex reductive group, let $B \subseteq G$ be a Borel subgroup, and let $T \subseteq B$ be a maximal torus. 
As mentioned above, a homogeneous space $X_0=G/H$ is called \emph{spherical} if $B$ acts on $X_0$ with a dense open orbit; we also say that $H$ is a \emph{spherical} subgroup of $G$ in this case. Conjugacy classes of spherical subgroups accept a combinatorial description in terms of discrete invariants, known as homogeneous spherical data, that were introduced by Luna (see \cite[Section 2]{Lun01} or \cite[Section 30.11]{Tim11}) under the names of \emph{homogeneous spherical data}. Roughly speaking, it is a quadruple $(\Lambda,\Pi^p,\Sigma,\DDD^a)$ with $\Lambda \subseteq \X$ a sublattice, $\Pi^p$ and $\Sigma$ two finite subsets of $\X$, and $\DDD^a$ a finite multisubset of the dual of the sublattice  of $\X$ generated by $\Sigma$, which satisfy some compatibility conditions (see \cite[Definition 30.21]{Tim11}). 
Let us also mention that, for any spherical subgroup $H \subseteq G$, the quotient group $N_G(H)/H$ is a diagonalizable group (\cite[Section 5.2]{BP87}); in particular, it is an abelian group. 

\smallskip

Let us fix a real group structure $\sigma$ on $G$. The $\star$-action induced by $\sigma$ on the based root datum associated to the triple $(G,B,T)$ (see Definition \ref{def: star-action}) induces a $\Gamma$-action on the set of homogeneous spherical data (see \cite[Sections 2.15-2.16]{BG18} for details, but beware that the notation is different from \cite[Section 30.11]{Tim11}).

The next result gives a complete answer to Question \ref{item: existence part of Key-Prob3} of the Key Problem for spherical homogeneous spaces.

\begin{theorem}\emph{(\cite[Main Theorem 1.13]{BG18} with $k_0=\R$)}\label{th: main th BG}\\
We keep the previous notation.
Let $X_0=G/H$ be a spherical homogeneous space. 
Then $X_0$ admits a $(G,\sigma)$-equivariant real structure if and only if 
\begin{enumerate}
\item\label{item: th 4.1 (i)} the homogeneous spherical datum of the conjugacy class of $H$ is preserved by the $\star$-action induced by $\sigma$; and
\item\label{item: th 4.1 (ii)} a certain cohomology class $\kappa_*(\delta([\sigma])) \in \H^2(\Gamma,N_G(H)/H)$ is trivial (see \cite[Section 3]{BG18} for the precise definition of this cohomology class).
\end{enumerate}
\end{theorem}

\begin{remark}\ 
\begin{itemize}
\item Condition \ref{item: th 4.1 (i)} of Theorem \ref{th: main th BG} is equivalent to the fact that $\sigma(H)$ is conjugate to $H$ in $G$, which is condition \ref{eq: sigma compatible} of Lemma \ref{lem: two conditions}.
\item If $Z(G) \subseteq H$ or if $\sigma$ is \emph{quasi-split} (i.e.~stabilizes a Borel subgroup $B \subseteq G$), then the cohomology class $\kappa_*(\delta([\sigma]))$ is trivial, and so $X_0$ admits a $(G,\sigma)$-equivariant real structure if and only if condition \ref{item: th 4.1 (i)} of Theorem \ref{th: main th BG} holds. Recall that every real group structure on $G$ admits an inner twist that stabilizes the Borel subgroup $B \subseteq G$.
\item We observe that Theorem \ref{th: main th BG} looks very similar to Theorem \ref{th: real structures on linear reps}, but we do not know what a common generalization of these two results could be.
\end{itemize}
\end{remark}

\begin{example}{(\cite[Example 2.2]{MJT19})}\label{example A}
Let $G=\SL_n^{\times 3}$ with $n \geq 2$, and let $\sigma\colon (g_1,g_2,g_3)\mapsto (\overline{g_2}, \overline{g_1},{\overline{\leftexp{t}g_3^{-1}}})$, which can easily verified to be a real group structure on $G$.
Let 
\[H:=\{(h_1,h_1,h_2);\ h_1 \in \SL_n \text{ and } h_2 \in \mathrm{SO}_n\}\  \text{and}\  
H':=\psi(H),\] 
where $\psi\colon G \to G,\ (g_1,g_2,g_3) \mapsto (g_3,g_2,g_1)$ is an outer automorphism of $G$.
Then $H$ and $H'$ are both spherical subgroups of $G$, but the $\star$-action induced by $\sigma$ stabilizes only the homogeneous spherical datum of the conjugacy class of $H$.
Indeed, denoting by $\{\alpha_1^i,\ldots, \alpha_{n-1}^i\}$ and $\{\varpi_1^i,\ldots,\varpi_{n-1}^i\}$ the sets of simple roots and fundamental weights of the $i$-th factor $\SL_n$ ($i \in \{1,2,3\}$) of $G=\SL_n^{\times 3}$, the homogeneous spherical datum of the conjugacy class of $H$ is given by
\[
 \left\{
    \begin{array}{lll}
       \Lambda&:= &\Z\langle \varpi_1^1+\varpi_1^2,\ldots, \varpi_{n-1}^1+\varpi_{n-1}^2, \alpha_1^3, \ldots, \alpha_{n-1}^3 \rangle \\
        \Pi^p&:= &\varnothing\\
        \Sigma&:=&\{ \alpha_1^1+\alpha_1^2, \ldots, \alpha_{n-1}^1+\alpha_{n-1}^2, \alpha_1^3, \ldots, \alpha_{n-1}^3 \}\\
        \DDD^a&:=&\{ D_{1,-}^3, D_{1,+}^3, D_{2,-}^3, D_{2,+}^3,\ldots, D_{n-1,-}^3, D_{n-1,+}^3\}
    \end{array},
\right.
\]
while the homogeneous spherical datum of the conjugacy class of $H'$ is given by
\[
 \left\{
    \begin{array}{lll}
       \Lambda&:= &\Z\langle \varpi_1^2+\varpi_1^3,\ldots, \varpi_{n-1}^2+\varpi_{n-1}^3, \alpha_1^1, \ldots, \alpha_{n-1}^1 \rangle \\
        \Pi^p&:= &\varnothing\\
        \Sigma&:=&\{ \alpha_1^2+\alpha_1^3, \ldots, \alpha_{n-1}^2+\alpha_{n-1}^3, \alpha_1^1, \ldots, \alpha_{n-1}^1 \}\\
        \DDD^a&:=&\{ D_{1,-}^1, D_{1,+}^1, D_{2,-}^1, D_{2,+}^1,\ldots, D_{n-1,-}^1, D_{n-1,+}^1\}
    \end{array}.
\right.
\]
Moreover, the cohomology class $\kappa_*(\delta([\sigma])) \in \H^2(\Gamma,N_G(H)/H)$ is trivial, and so $X_0=G/H$ admits a $(G,\sigma)$-equivariant real structure by Theorem \ref{th: main th BG}.
\end{example}

\begin{example}{(\cite[Example 11.5]{BG18})}\label{example B}
Let $G=\SO_{10}$, let $\sigma$ be any real group structure on $G$, and let $H=\SO_9$ embedded in $G$ in the standard way. Then $X_0=G/H$ is a spherical homogeneous space, and the conjugacy class of $H$ is stable for the $\star$-action induced by $\sigma$.
Indeed, denoting by $\{\alpha_1,\ldots, \alpha_5\}$ and $\{\varpi_1,\ldots,\varpi_5\}$ the sets of simple roots and fundamental weights of $G=\SO_{10}$ (with Bourbaki's notation), the homogeneous spherical datum of the conjugacy class of $H$ is given by
\[
 \left\{
    \begin{array}{lll}
       \Lambda&:= &\Z\langle \varpi_1 \rangle \\
        \Pi^p&:= & \{\alpha_2,\alpha_3,\alpha_4,\alpha_5\}\\
        \Sigma&:=&\{ \varpi_1\}\\
        \DDD^a&:=&\varnothing
    \end{array},
\right.
\]
and the Galois group acts on the character lattice $\X=\Z\langle \varpi_1,\ldots,\varpi_5\rangle$ either trivially or by multiplying $\varpi_5$ by $-1$ (which corresponds to swapping $\alpha_4$ and $\alpha_5$).

Let $\sigma'$ be the real group structure on $G$ whose real locus is the non-compact real Lie group usually denoted by $\SO^*(10)$. If $\sigma$ belongs to the equivalence class of $\sigma'$, then the cohomology class $\kappa_*(\delta([\sigma])) \in \H^2(\Gamma,N_G(H)/H)$ is nontrivial, and so $X_0$ does not admit a $(G,\sigma)$-equivariant real structure. 
For any other $\sigma$, the cohomology class $\kappa_*(\delta([\sigma]))$ is trivial, and so $X_0$ admits a $(G,\sigma)$-equivariant real structure. 
\end{example}

\smallskip

Let us now briefly consider Question \ref{item: quantity part of Key-Prob3} of the Key Problem for spherical homogeneous spaces.
It follows from Corollary \ref{cor: finiteness for homog spaces} that, if the spherical homogeneous space $X_0=G/H$ admits a $(G,\sigma)$-equivariant real structure, then $X_0$ admits a finite number of them, up to equivalence.  Moreover, as mentioned earlier, the group $\Aut^G(X_0) \simeq N_G(H)/H$ is diagonalizable, and so it fits into an exact sequence of abelian $\Gamma$-groups
\[
0 \to \Aut^G(X_0)^\circ \to \Aut^G(X_0) \to \pi_0(\Aut^G(X_0)) \to 0,
\]
where $\Aut^G(X_0)^\circ$  is a torus and $\pi_0(\Aut^G(X_0))$ is a finite abelian group, which makes it possible in practice to compute $\H^1(\Gamma,\Aut^G(X_0))$ through a long exact sequence in Galois cohomology. (See \cite[Section 3.4]{MJT18} and \cite[Section 3]{MJT19} for some explicit examples of $\H^1(\Gamma,\Aut^G(X_0))$ computations.)

\begin{works}\ 
\smallskip
\begin{itemize}
\item In \cite{BG18}, Borovoi and Gagliardi study more generally descent data and $\K$-forms for spherical varieties, with $\K$ an arbitrary base field of characteristic zero (but not specifically the field of real numbers).
\smallskip
\item Theorem \ref{th: main th BG} was first proved for \emph{horospherical} homogeneous spaces (i.e. homogeneous spaces $G/H$ with $H$ containing a maximal unipotent subgroup of $G$) by Moser-Jauslin--Terpereau in \cite{MJT18}, and then later generalized by Borovoi--Gagliardi to arbitrary spherical homogeneous spaces over a base field of characteristic zero.\smallskip
\item Equivariant real structures on spherical homogeneous spaces $X_0=G/H$, under the extra assumption that $\Aut^G(X_0) \simeq N_G(H)/H$ is finite, have been studied by Akhiezer in \cite{Akh15}, by Akhiezer--Cupit-Foutou in \cite{ACF14}, by Cupit-Foutou in \cite{CF15}, by Borovoi in \cite{Bor20}, and by Snegirov in \cite{Sne20}.
\smallskip 
\item In the particular case where $X_0=G/H$ is a \emph{symmetric space} (i.e.~$G^\theta \subseteq H \subseteq N_G(G^\theta)$ with $\theta \in \Aut_{\gr}(G)$ a group involution), a practical criterion for the existence of an equivariant real structure on $X_0$, using the involution $\theta$ instead of the homogeneous spherical datum of the conjugacy class of $H$, has been obtained by Moser-Jauslin--Terpereau in \cite{MJT19} .
\smallskip
\item Let $X_0=G/H$ be a spherical homogeneous space with a $(G,\sigma)$-equivariant real structure $\mu$ such that $X_0(\C)^\mu$ is non-empty. Then the real Lie group $G(\C)^\sigma$ acts on $X_0(\C)^\mu$ with finitely many orbits. 
When $X_0$ is a symmetric space, a combinatorial description of these orbits using Galois cohomology has been obtained by Cupit-Foutou--Timashev in \cite{CFT18} (see also \cite[Chp.~6]{BJ06}), and when $X_0$ is arbitrary but $\sigma$ is split, a parametrization of these orbits through geometric methods has been obtained by Cupit-Foutou--Timashev in \cite{CFT19}.
\end{itemize}
\end{works}

\subsection{Extension of real structures to spherical embeddings}\label{subsec: extension of real structures to sph embeddings}
In the previous subsection, we investigated the Key Problem for spherical homogeneous spaces. In this subsection, we go to the next step and consider the Key Problem for spherical embeddings (although, historically, things have rather worked the other way round). More precisely, given a spherical homogeneous space $X_0=G/H$ and a $G$-equivariant embedding $X_0 \hookrightarrow X$, we give the combinatorial criterion obtained by Huruguen in \cite{Hur11} to determine when an equivariant real structure on $X_0$ extends to $X$, generalizing what was done in Section  \ref{subsec: real structures on toric varieties} for toric varieties. For brevity, we do not give details regarding the theory of spherical embeddings; the interested reader is referred to \cite{Kno91,Tim11,Per14} for a detailed presentation.

\smallskip

We keep the same notation as in the previous subsection, and we fix a spherical homogeneous space $X_0=G/H$. 
The \emph{weight lattice} of $X_0$ is
\[ \X(X_0)=\{ \chi \in \X(B) \ | \ \C(X_0)_\chi^{(B)} \neq \{0\} \} \ \subseteq \X(B):=\Hom_{\gr}(B,\G_m)\simeq \X(T), \]
where 
\[\C(X_0)_\chi^{(B)}=\{ f \in \C(X_0)\ |\ \forall b \in B,\ b \cdot f= \chi(b)f\}\ \text{ for }\ \chi \in \X(B).\] 
It is a free abelian group of finite rank. 
The (finite) set of \emph{colors} of $X_0$ is 
\[\DDD(X_0)=\{ \text{$B$-stable prime divisors of $X_0$} \}.\]
A \emph{colored cone} for $X_0$ is a pair formed by a strictly convex polyhedral cone in $\X(X_0)_\Q^\vee$ and a subset of $\DDD(X_0)$ satisfying some conditions.
A \emph{colored fan} for $X_0$ is a finite collection of colored cones for $X_0$ satisfying some extra conditions. (See \cite[Section 3]{Kno91} for details.)

\begin{theorem}\emph{(\cite[Theorem 3.3]{Kno91})}
There is a bijection between isomorphism classes of $G$-equivariant embeddings of $X_0$ and colored fans for $X_0$.
\end{theorem}

Let now $\sigma$ be a real group structure on $G$, and let $\mu$ be a $(G,\sigma)$-equivariant real structure on $X_0$. Then $\mu$ induces a $\Gamma$-action on the weight lattice $\X(X_0)$ and on the set of colors $\DDD(X_0)$; see \cite[Section 2.2]{Hur11} for a detailed description of this $\Gamma$-action. 
This $\Gamma$-action on  $\X(X_0)$ and $\DDD(X_0)$ in turn induces a $\Gamma$-action on the set of colored fans for $X_0$. 

\begin{remark}
A somewhat surprising fact is that this $\Gamma$-action on the set of colored fans for $X_0$ does not depend on the choice of $\mu$, only on $\sigma$, and can be recovered from the $\star$-action induced by $\sigma$ on the based root datum associated to the triple $(G,B,T)$; see \cite[Section 7]{BG18} for details. 
\end{remark}

\begin{theorem}\label{th:real structure extend to spherical embeddings} \emph{(\cite[Theorem 2.23]{Hur11}, see also \cite[Theorem 9.1]{Wed18})}\\
We keep the previous notation. Let $X_0 \hookrightarrow X$ be a $G$-equivariant embedding of $X_0$.
A $(G,\sigma)$-equivariant real structure $\mu$ on $X_0=G/H$ extends to $X$ if and only if the corresponding colored fan is $\Gamma$-stable.
\end{theorem}

\begin{remark}\
\begin{itemize}
\item Theorem \ref{th:real structure extend to toric embeddings} for toric embeddings is a particular case of Theorem \ref{th:real structure extend to spherical embeddings}.
\item Unlike the toric case, the Galois descent $\C/\R$ is not always effective for spherical varieties. This means that a $(G,\sigma)$-equivariant real structure on a spherical variety $X$ does not always correspond to a real $(G/\langle \sigma \rangle)$-form of $X$. We refer to \cite[Section 2.4]{Hur11} for an explicit example where this situation occurs. The property that $X$ admits a covering by $\Gamma$-stable quasi-projective open subsets can however be expressed combinatorially (see \cite[Corollary 3.2.12]{Per14}). See also \cite[Proposition 2.27]{Hur11} for a list of sufficient conditions to guaranty effectiveness.
\end{itemize}
\end{remark}

\begin{example}(\cite[Example 3.34]{MJT18})
Let $U$ be a maximal unipotent subgroup of $\SL_2$. Then $X_0=\SL_2/U \simeq \A_\C^2 \setminus \{0\}$ is a spherical homogeneous space whose nontrivial equivariant embeddings and corresponding colored fans are given as follows:\ 
$\A_\C^2$, \ $\P_\C^2 \setminus \{0\}$, \ $\Bl_0(\A_\C^2)$, \ $\P_\C^2$,  and $\Bl_0(\P_\C^2)$.
\vspace{3mm}
  \begin{center}
\includegraphics[width=10cm]{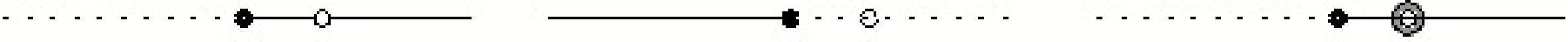}\vspace{2mm}
\includegraphics[width=10cm]{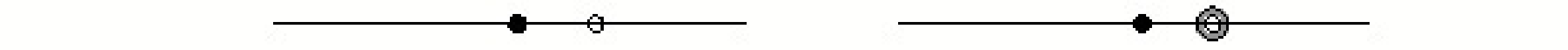}
\end{center}
\vspace{2mm}
(The small white dot corresponds to the unique color of $\SL_2/U$.) 
We saw in Example \ref{ex:real group structures on SL2} that there are two inequivalent real group structures on $\SL_2$, which we denoted by $\sigma_s$ and $\sigma_c$. 
It follows from Theorem \ref{th: main th BG} and Proposition \ref{prop:Galois H1 to param eq real structures} that there exists a unique equivalence class of $(\SL_2,\sigma_s)$-equivariant real structures on $X_0$, but that there is no $(\SL_2,\sigma_c)$-equivariant structure on $X_0$ as the cohomology class $\kappa_*(\delta([\sigma_c])) \in \H^2(\Gamma,N_G(U)/U)$ is nontrivial.
Moreover, the $\star$-action induced by $\sigma_s$ is trivial (see Remark \ref{rk: star action trivial in the split case}), and so any $(\SL_2,\sigma_s)$-equivariant real structure on $X_0$ extends to any $\SL_2$-equivariant embedding $X_0 \hookrightarrow X$. All $(\SL_2,\sigma_s)$-equivariant real structures on $X$ are furthermore equivalent since $\Aut^G(X) \simeq \Aut^G(X_0) (\simeq \G_{m,\C})$.
\end{example}

\smallskip

Let us now move on to Question \ref{item: quantity part of Key-Prob3} of the Key Problem for a given spherical variety $X$ with open orbit $X_0$. As for toric varieties, once the real group structure $\sigma$ on $G$ is fixed, we have the following alternative: Either all $(G,\sigma)$-equivariant real structures on $X_0$ extend to $X$, or none of them do.  
However, contrary to the toric case, it may happen that two equivalent $(G,\sigma)$-equivariant real structures on $X_0$ extend to inequivalent $(G,\sigma)$-equivariant real structures on $X$; indeed, the group inclusion $\Aut^G(X) \hookrightarrow \Aut^G(X_0)$ may not be an isomorphism. 
Consequently, to answer Question \ref{item: quantity part of Key-Prob3} of the Key Problem, the best option seems again to write a long exact sequence in Galois cohomology as in the homogeneous case.

\newpage

\begin{works}\ 
\begin{itemize}
\smallskip
\item In \cite{Hur11}, Huruguen works over any perfect base field, not necessarily the field of real numbers.  Let us note that Huruguen assumes that the equivariant real structure $\mu$ on $X_0$ satisfies $X_0(\C)^\mu \neq \varnothing$, but this condition plays no role in the proof of Theorem \ref{th:real structure extend to spherical embeddings}.
The results of Huruguen have been later generalized by Wedhorn in \cite{Wed18} to the case where the base field is arbitrary. Also, Wedhorn prefers to work in the category of algebraic spaces instead of schemes, which allows him to avoid the problem of effectiveness for the Galois descent.

\smallskip

\item Akhiezer and Cupit-Foutou study in \cite{Akh15,ACF14,CF15} the extension problem for an equivariant real structure on a spherical homogeneous space $X_0$ admitting a wonderful compactification $X_0 \hookrightarrow X$.

\smallskip

\item In addition to proving Theorem \ref{th: main th BG}, Borovoi--Gagliardi prove in \cite[Section 7]{BG18} a reformulation of part of the results obtained by Huruguen in \cite{Hur11}. They also provide many examples in \cite[Section 11]{BG18} of spherical varieties for which they apply Theorem \ref{th: main th BG} to answer Question \ref{item: existence part of Key-Prob3} of the Key Problem.

\smallskip

\item As a byproduct of Theorems \ref{th: main th BG} and \ref{th:real structure extend to spherical embeddings}, a classification of  the equivariant real structures on smooth projective horospherical varieties of Picard rank $1$ (classified by Pasquier in \cite{Pas09}) has been obtained by Moser-Jauslin--Terpereau in \cite[Section 3.6]{MJT18}. 

\end{itemize}
\end{works}

\section{Real structures on almost homogeneous \texorpdfstring{$\SL_2$}{SL2}-threefolds} \label{sec:almost homg SL2-threefolds}
Let $G$ be a complex reductive group, and let $B$ be a Borel subgroup of $G$. We recall that the \emph{complexity} of a $G$-variety is the codimension of a general $B$-orbit. For instance, the complexity-zero varieties are precisely the spherical varieties considered in Section \ref{sec: real structures on spherical varieties}. In this section, based on the work of Moser-Jauslin--Terpereau \cite[Section 4]{MJT20}, we study the Key Problem for a family of complexity-one varieties admitting a combinatorial description quite similar to that of spherical embeddings, namely the \emph{almost homogeneous $\SL_2$-threefolds}. Examples include the rank $1$ Fano threefolds $\P_\C^3$, $Q_3$, $X_5$ and $X_{22}^{\textrm{MU}}$ (with the notation of \cite[Section 12.2]{IP99}), and the $\P_\C^1$-bundle $\P(T_{\P_\C^2})$.

\smallskip

We have seen in Example \ref{ex:real group structures on SL2} that any real group structure on $\SL_2$ is equivalent to either $\sigma_s\colon\  g \mapsto \overline{g}$ or $\sigma_c\colon\  g \mapsto {}^{t}\overline{g^{-1}}$, with corresponding real loci $\SL_2(\R)$ and $\SU_2(\C)$ respectively. 
Moreover, $\SL_2$ has no outer automorphism, and so we can assume without loss of generality that $\sigma \in \{ \sigma_s,\sigma_c\}$ when studying $(\SL_2,\sigma)$-equivariant real structures on $\SL_2$-varieties (see the second part of \cite[Proposition 1.1]{MJT20}). Let us also note that, for any real group structure $\sigma$ on $\SL_2$, the $\star$-action induced by $\sigma$ is trivial (see Remark \ref{rk: star action trivial in the split case}), and so the $\star$-action will play no role when considering the Key Problem for almost homogeneous $\SL_2$-threefolds (contrary to the spherical case studied in Section \ref{sec: real structures on spherical varieties}, where the $\star$-action was playing a leading role to answer Question \ref{item: existence part of Key-Prob3} of the Key Problem). 

Let $X$ be an almost homogeneous $\SL_2$-threefold. Then it contains a dense open orbit $X_0$ isomorphic to $\SL_2/H$ with $H \subseteq \SL_2$ a finite subgroup; these are well-known (see \cite{Klein93}): there are the cyclic groups of order $n$ (conjugate to $A_n$), the binary dihedral groups of order $4n-8$ (conjugate to $D_n$ with $n \geq 4$), and the binary polyhedral groups (conjugate to $E_n$ with $n \in \{6,7,8\}$).

\begin{theorem}\emph{(\cite[Theorem~1.6]{MJT20})}\label{th:D}
Let $H$ be a finite subgroup of $\SL_2$, and let $\sigma$ be a real group structure on $\SL_2$.
Then $X_0=\SL_2/H$ admits an $(\SL_2,\sigma)$-equivariant real structure. Moreover, the equivalence classes of $(\SL_2,\sigma)$-equivariant real structures on $X_0$ and their real loci are listed in \cite[Appendix C]{MJT20}; there are $1$, $2$ or $3$ equivalence classes depending on $H$ and $\sigma$.
\end{theorem}

Following the strategy outlined in Section \ref{sec: equivariant real structures} to describe the equivariant real structures on almost homogeneous varieties, we now consider the problem of determining when an $(\SL_2,\sigma)$-equivariant real structure on $X_0=\SL_2/H$ extends to a given $\SL_2$-equivariant embedding $X_0 \hookrightarrow X$. 
As for spherical embeddings, there exists a combinatorial classification of the $\SL_2$-equivariant embeddings of $X_0$. (In fact, both descriptions are particular cases of the Luna-Vust theory to classify equivariant embeddings of arbitrary homogeneous spaces.) Consider the sets
\[
\V^{\SL_2}=\V^{\SL_2}(X_0)=\{\text{$\SL_2$-invariant geometric valuations\footnotemark \ of $\C(X_0)$}\}
\]
and, for a given Borel subgroup $B$ of $\SL_2$,
\footnotetext{ Recall that a (discrete) \emph{valuation} of $\C(X_0)$ is a group homomorphism $\nu\colon (\C(X_0)^*,\times) \to (\Q,+)$ satisfying $\nu(a+b) \geq \min(\nu(a),\nu(b))$ when $a+b \neq 0$, whose kernel contains $\C^*$, and whose image is a discrete subgroup of $(\Q,+)$. A valuation $\nu$ is called \emph{geometric} if $\nu=c\cdot \nu_D$ for some $c \in \Q_+^*$ and some divisor $D$ on a $G$-variety equivariantly birational to $X_0$, where $\nu_D(f)$ is the order of vanishing of $f$ along $D$.}
\[
\DDD^B=\DDD^B(X_0)=\{ \text{$B$-stable prime divisors of $X_0$} \} \simeq B \backslash \SL_2/H\simeq \P_\C^1/H.
\]
The idea of the classification of the $\SL_2$-equivariant embeddings of $X_0$ is that any such embedding $X_0 \hookrightarrow X$ corresponds to a collection of \emph{colored data}, i.e.~a collection of pairs $(\W_i,\RR_i)_{i \in I}$, with $\W_i \subseteq \V^{\SL_2}$ and $\RR_i \subseteq \DDD^B$, satisfying some conditions.
In fact, each colored datum corresponds to an $\SL_2$-orbit of $X$, and the conditions that the collection of pairs $(\W_i,\RR_i)_{i \in I}$ must satisfy encode the fact that the union of the corresponding $\SL_2$-orbits is indeed an $\SL_2$-variety. Also, $\SL_2$-orbits are divided into six different types (depending on their colored data) denoted by 
\begin{itemize}
\item $\AA$, $\AB$, $\BB_{+}$, or $\BB_{-}$ for orbits isomorphic to $\P_\C^1$;
\item $\BB_0$ for fixed points; and
\item $\CC$ for $2$-dimensional orbits.
\end{itemize}   
This combinatorial classification was first given by Luna--Vust in \cite[Section 9]{LV83} for equivariant embeddings of $\SL_2$, and then later generalized to equivariant embeddings of $\SL_2/H$, with $H \subseteq \SL_2$ a finite subgroup, by Moser-Jauslin in \cite{MJ87,MJ90} and Bousquet in \cite{Bou00}; we refer to \cite[Appendix B]{MJT20} for a self-contained summary. 

\begin{remark}Let us mention that there is an alternative description of equivariant embeddings of $\SL_2/H$, in terms of \emph{colored hypercones} and \emph{colored hyperfans}, due to Timashev; see \cite[Section 16]{Tim11}.
\end{remark}

In practice, the collection of colored data corresponding to a given equivariant embedding of $\SL_2/H$ is represented by a \emph{skeleton diagram}. For instance the skeleton diagrams corresponding to the varieties $Q_3$, $V_5$ and $V_{22}^{\textrm{MU}}$ are respectively

\begin{figure}[h!]
\begin{tikzpicture}[scale=1.2]
\path (0,0) coordinate (origin);
\draw (0,0) circle (2pt) ;
\path (0:1.7cm) coordinate (P0);
\path (0:1.1cm) coordinate (P0););
\path (75:1.3cm) coordinate (P1););
\path (19:.75cm) coordinate (P2);
\path (310:1.3cm) coordinate (P3);
\path (55:.75cm) coordinate (P4);
\path (-30:.75cm) coordinate (P5);
\path (90:.8cm) coordinate (Q0);
\draw[line width=0.3mm]  (origin) -- (P0) (origin) -- (P3) (origin) -- (P4) (origin) -- (P5) (origin) -- (P1) (origin) -- (P2) ;
\fill[white] (0,0) circle (1.5pt) ;
\node at (.5,.5){$\vdots$};
\node at (.6,.8){\tiny{$-\frac{5}{6}$}};
\node at (0,1.3){\tiny{$-\frac{1}{2}$}};
\node at (.4,-1.1){\tiny{$-\frac{1}{2}$}};
\node at (1.5,0){\tiny{$-\frac{2}{3}$}};
\node at (-.2,-.2){\tiny{$-1$}};
\node at (.3,-2){\tiny{$H=E_6$}};
\draw[line width=0.3mm] (P1) -- +(165:3pt)  (P1)-- +(345:3pt); 
\end{tikzpicture}
\begin{tikzpicture}[scale=1.2]
\path (0,0) coordinate (origin);
\draw (0,0) circle (2pt) ;
\path (0:1.7cm) coordinate (P0);
\path (0:1.1cm) coordinate (P0););
\path (75:1.3cm) coordinate (P1););
\path (19:.75cm) coordinate (P2);
\path (310:1.3cm) coordinate (P3);
\path (55:.75cm) coordinate (P4);
\path (-30:.75cm) coordinate (P5);
\path (90:.8cm) coordinate (Q0);
\draw[line width=0.3mm] (origin) -- (P0) (origin) -- (P3) (origin) -- (P4) (origin) -- (P5) (origin) -- (P1) (origin) -- (P2) ;
\fill[white] (0,0) circle (1.5pt) ;
\node at (.5,.5){$\vdots$};
\node at (.6,.8){\tiny{$-\frac{5}{6}$}};
\node at (0,1.3){\tiny{$-\frac{1}{2}$}};
\node at (.4,-1.1){\tiny{$-\frac{1}{2}$}};
\node at (1.5,0){\tiny{$-\frac{2}{3}$}};
\node at (-.2,-.2){\tiny{$-1$}};

\node at (.3,-2){\tiny{$H=E_6$}};
\draw[line width=0.3mm] (P3) -- +(40:3pt)  (P3)-- +(220:3pt); 
\end{tikzpicture}
\begin{tikzpicture}[scale=1.2]
\path (0,0) coordinate (origin);
\draw (0,0) circle (2pt) ;
\path (0:1.7cm) coordinate (P0);
\path (0:1.1cm) coordinate (P0););
\path (75:1.1cm) coordinate (P1););
\path (19:.75cm) coordinate (P2);
\path (310:1.5cm) coordinate (P3);
\path (55:.75cm) coordinate (P4);
\path (-30:.75cm) coordinate (P5);
\path (90:.8cm) coordinate (Q0);
\draw[line width=0.3mm] (origin) -- (P0) (origin) -- (P3) (origin) -- (P4) (origin) -- (P5) (origin) -- (P1) (origin) -- (P2) ;
\fill[white] (0,0) circle (1.5pt) ;
\node at (.5,.5){$\vdots$};
\node at (.6,.8){\tiny{$-\frac{11}{12}$}};
\node at (0,1.3){\tiny{$-\frac{5}{6}$}};
\node at (.4,-1.1){\tiny{$-\frac{2}{3}$}};
\node at (1.5,0){\tiny{$-\frac{3}{4}$}};
\node at (-.2,-.2){\tiny{$-1$}};

\node at (.3,-2){\tiny{$H=E_7$}};
\draw[line width=0.3mm] (P3) -- +(40:3pt)  (P3)-- +(220:3pt); 
\end{tikzpicture}
\end{figure}

The choice of a pair $(\sigma,\mu)$, with $\sigma$ a real group structure on $\SL_2$ and $\mu$ an $(\SL_2,\sigma)$-equivariant real structure on $X_0$, induces a $\Gamma$-action on the sets $\V^{\SL_2}$ and $\DDD^B$ (see \cite[Section 4.3]{MJT20} for details). This $\Gamma$-action can be visualized on the skeleton diagrams; indeed, it corresponds to a permutation of certain spokes of the same length.
The following result provides a complete answer to Question \ref{item: existence part of Key-Prob3} of the Key Problem for almost homogeneous $\SL_2$-threefolds.

\begin{theorem}\emph{(\cite[Theorem~4.7]{MJT20})}\label{th:E} 
We keep the previous notation. Let $X_0 \hookrightarrow X$ be an $\SL_2$-equivariant embedding of $X_0=\SL_2/H$.
An $(\SL_2,\sigma)$-equivariant real structure $\mu$ on $X_0$ extends to an $(\SL_2,\sigma)$-equivariant real structure $\widetilde\mu$ on $X$ if and only if the  $\Gamma$-actions on $\V^{\SL_2}$  and $\DDD^B$ induced by the pair $(\sigma,\mu)$ stabilize the collection of colored data corresponding to the $\SL_2$-orbits of $X$. Moreover, the real structure $\widetilde\mu$ on $X$ is effective if and only if every $\SL_2$-orbit of $X$ of type $\BB_0$ or $\BB_-$ is fixed by the $\Gamma$-action.
\end{theorem}

\begin{remark}
It must be stressed that, contrary to the spherical case (see Section \ref{sec: real structures on spherical varieties}), the $\Gamma$-action on the colored equipment of $X_0=\SL_2/H$ depends not only on $\sigma$, but also on $\mu$. In fact, it is even possible, for a given $\sigma$, to have two equivalent $(\SL_2,\sigma)$-equivariant real structures on $X_0$ such that only one of them extends to a given $\SL_2$-equivariant embedding $X_0 \hookrightarrow X$.
\end{remark}

\begin{example}(\cite[Example~4.11]{MJT20})\label{ex:first bis example}
Let $X=\P_\C^1\times\P_\C^1 \times \P_\C^1$ on which $\SL_2$ acts diagonally. Then the stabilizer of the point $x=([1:1],[1:0],[0:1])$ is $H=\{ \pm I_2 \}=A_2$, and so $(X,x)$ is an $\SL_2$-equivariant embedding of $\SL_2/H=\PGL_2$. 
The orbit decomposition of $X$ is $\ell \sqcup S_1 \sqcup S_2 \sqcup S_3 \sqcup X_0$, where $X_0 \simeq \PGL_2$ is the dense open orbit, $S_i \simeq \P_\C^1 \times \P_\C^1 \setminus \Delta$, and $\ell \simeq \P_\C^1$.
Let us note that
\[
\mathfrak{S}_3 \simeq \Aut^{\SL_2}(X) \hookrightarrow \Aut^{\SL_2}(X_0) \simeq \PGL_2(\C), \ (12) \mapsto \begin{bmatrix}
 i & 0 \\i & -i
\end{bmatrix} \ \text{and}  \ (23) \mapsto \begin{bmatrix}
0 & i \\ i & 0
\end{bmatrix},
\]
where the symmetric group $\mathfrak{S}_3$ acts on $X$ by permuting the three factors. 

We fix a real group structure $\sigma$ on $\SL_2$. By \cite[Theorem 1.6]{MJT20}, there are exactly two equivalence classes of $(\SL_2,\sigma)$-equivariant real structures on $X_0$. 
Using Theorem \ref{th:E}, one can show that there always exists an $(\SL_2,\sigma)$-equivariant real structure on $X_0$ that extends to $X$, and using Proposition \ref{prop:Galois H1 to param eq real structures}, a direct computation of $\H^1(\Gamma,\Aut^{\SL_2}(X))$ yields that $X$ admits exactly two equivalence classes of $(\SL_2,\sigma)$-equivariant real structures.
Moreover, all $(\SL_2,\sigma)$-equivariant real structures on $X$ restrict to equivalent $(\SL_2,\sigma)$-equivariant real structures on $X_0$.
\end{example}

Other examples where Theorem \ref{th:E} has been applied to determine the equivariant real structures on certain almost homogeneous $\SL_2$-threefolds can be found in \cite[Section 4.3]{MJT20}. In particular, \cite[Example 4.8]{MJT20} gives an example of an almost homogeneous $\SL_2$-threefold with an equivariant real structure that is \emph{not} effective.
Also, a classification of the equivariant real structures on minimal smooth completions of $X_0=\SL_2/H$, when $H$ is non-cyclic, has been obtained in \cite[Section 4.4]{MJT20}. (Here, we refer to any $\SL_2$-equivariant embedding
	$X_0 \hookrightarrow X$ as a \emph{minimal smooth completion} of $X_0$, provided that $X$ is a smooth complete variety and any birational $\SL_2$-equivariant morphism $X \to X'$, with $X'$ smooth, is an isomorphism.)

\begin{works}\ 
\begin{itemize}
\smallskip 
\item Examples of smooth compact connected three-dimensional algebraic spaces that are not schemes and on which $\SL_2$ acts with a dense open orbit have been obtained by Luna--Moser-Jauslin--Vust in \cite{LMV89}. (One way to produce such examples in the real setting is to take the quotient $X/\langle \mu \rangle$ with $X$ a complete almost homogeneous $\SL_2$-threefold and $\mu$ an $(\SL_2,\sigma)$-equivariant real structure on $X$ that is not effective.)

\smallskip

\item Let $G$ be any complex reductive group, and let $\sigma$ be a real group structure on $G$. A combinatorial criterion to determine if a $(G,\sigma)$-equivariant real structure $\mu$ on an arbitrary homogeneous space $X_0=G/H$ extends to a given $G$-equivariant embedding $X_0 \hookrightarrow X$, relying on the classical Luna--Vust theory (see \cite{LV83}), has been obtained by Moser-Jauslin--Terpereau in \cite{MJT20}, generalizing Theorems \ref{th:real structure extend to spherical embeddings} and \ref{th:E}, but the combinatorics is then more involved.
\end{itemize}
\end{works}

\section{Some open questions}\label{sec: open questions}
We complete this survey paper with some open questions related to the Key Problem. Throughout this last section, we denote by $G$ a complex reductive group.

\smallskip

In Section \ref{sec:almost homg SL2-threefolds}, we have studied the Key Problem for almost homogeneous $\SL_2$-threefolds. On the other hand, there exists a combinatorial description for com\-plexity-one $G$-varieties due to Timashev, in terms of \emph{colored hypercones} and \emph{colored hyperfans} (see \cite[Section 16]{Tim11}), that generalizes the one given by Luna--Vust in \cite[Section 9]{LV83} for almost homogeneous $\SL_2$-threefolds.
\begin{question} \label{Q1}
Is it possible to extend the results obtained for almost homogeneous $\SL_2$-threefolds (see Section \ref{sec:almost homg SL2-threefolds}) to arbitrary complexity-one $G$-varieties?
\end{question}
(Let us mention that for affine varieties endowed with a complexity-one torus action, a positive answer to Question \ref{Q1} follows from the work of Langlois \cite{Lan15} and Gillard \cite{Gil20}; see Section \ref{subsec: real structures on T-varieties}).

\smallskip

In \cite[Section 4.4]{MJT20}, we have determined the equivariant real structures on the minimal smooth completions of $\SL_2/H$ when $H$ is a non-cyclic finite subgroup of $\SL_2$ (in which case the underlying $\SL_2$-variety is projective, and so any real structure is effective). But ``most" of the almost homogeneous $\SL_2$-threefolds that appear in the literature are actually related to minimal smooth completions of $\SL_2/H$ with $H$ cyclic.

\begin{question}  \label{Q2}
What are the effective equivariant real structures on the minimal smooth completions of $X_0=\SL_2/H$ when $H$ is cyclic?
\end{question}

It is certainly possible to try to answer Question \ref{Q2} using the same techniques as in \cite[Section 4.4]{MJT20} when $H$ is a non-cyclic finite subgroup of $\SL_2$. 
However, when $H$ is cyclic, the group $\Aut^{\SL_2}(X_0)$ is infinite, which makes it more complicated to compute $\H^1(\Gamma,\Aut^{\SL_2}(X_0))$, and there are between seven and eleven minimal smooth completions of $X_0$ to consider (depending on the cardinality of $H$). Moreover, the underlying $\SL_2$-varieties are not all projective, which makes the question of the effectiveness of the Galois descent nontrivial in this case.

\smallskip

Among all the complex Fano threefolds (see \cite[Section 12.2]{IP99} for the list of the smooth ones), many of them are varieties of complexity $\leq 1$. 
Therefore, results obtained in recent years concerning the equivariant real structures on spherical and complexity-one varieties should allow a complete classification of the real structures on these varieties.
\begin{question}
What are the (equivariant) real structures on the complex Fano threefolds of complexity $\leq 1$? And what are the corresponding real loci?
\end{question}

Let us now consider $X_0=G/H$ an arbitrary homogeneous space.
If $H$ is connected and $\dim(X_0) \leq 10$ or if $X_0$ is spherical, then it is known that $X_0$ is a rational variety (see \cite[Theorem 5.9]{CZ17} and \cite[Corollary 2.1.3]{Per14}). However, a real form of a rational variety is not necessarily rational.

\begin{question}
Given a rational homogeneous space $X_0=G/H$, what are the rational real forms of $X_0$? And is it possible to characterize rationality via the real locus? 
\end{question}

Let us now focus on the second part of Question \ref{item: quantity part of Key-Prob3} of the Key Problem. The automorphism group $\Aut^G(X)$ of a spherical $G$-variety $X$ is a diagonalizable group; in particular, it is a linear algebraic group, and so the set $\H^1(\Gamma,\Aut^G(X))$ is finite. On the other hand, there are examples of complexity-two varieties with uncountably many equivalence classes of equivariant real structures (see \cite{Les18,DO19,DFMJ21,DOY}). And in the middle lies the complexity-one case.

\begin{question}
Does a complexity-one variety always admit a finite number of equivalence classes of equivariant real structures?
\end{question}

Finally, for a random projective variety $X$, the identity component of its automorphism group $\Aut^\circ(X)$ is usually neither a linear algebraic group nor an abelian variety. It is therefore natural to consider the Key Problem without the assumption of reductivity, and even of linearity. 

\begin{question}
Let $K$ be a connected algebraic group (not necessarily linear), let $X$ be a complex $K$-variety, and let $\sigma$ be a real group structure on $K$.
Assume for instance that $X$ has general $K$-orbits of codimension $\leq 1$.
What can be said about the $(K,\sigma)$-equivariant real structures on $X$?
\end{question}

\bibliographystyle{alpha}
\bibliography{biblio}

\end{document}